# ASYMPTOTIC INFERENCE FOR HIGH-DIMENSIONAL DATA

By Jim Kuelbs and Anand N. Vidyashankar[1]

*University of Wisconsin and Cornell University*

In this paper, we study inference for high-dimensional data characterized by small sample sizes relative to the dimension of the data. In particular, we provide an infinite-dimensional framework to study statistical models that involve situations in which (i) the number of parameters increase with the sample size (that is, allowed to be random) and (ii) there is a possibility of missing data. Under a variety of tail conditions on the components of the data, we provide precise conditions for the joint consistency of the estimators of the mean. In the process, we clarify and improve some of the recent consistency results that appeared in the literature. An important aspect of the work presented is the development of asymptotic normality results for these models. As a consequence, we construct different test statistics for one-sample and two-sample problems concerning the mean vector and obtain their asymptotic distributions as a corollary of the infinite-dimensional results. Finally, we use these theoretical results to develop an asymptotically justifiable methodology for data analyses. Simulation results presented here describe situations where the methodology can be successfully applied. They also evaluate its robustness under a variety of conditions, some of which are substantially different from the technical conditions. Comparisons to other methods used in the literature are provided. Analyses of real-life data is also included.

**1. Introduction.** Modern scientific technology is providing a class of statistical problems that typically involve data that are high dimensional, and frequently lead to questions involving simultaneous inference for large sets of parameters. The number of parameters in these datasets is often random,

Received October 2008; revised February 2009.

[1]Supported in part by NSF Grant DMS-000-03-07057 and also by grants from the NDC Health Corporation.

and grows rapidly in comparison to the sample size; furthermore, there can be missing observations. Microarrays epitomize this situation, but similar problems arise in other areas such as polymerase chain reactions, proteomics, functional magnetic resonance imaging, and astronomy. For example, in microarray experiments the number of expressed genes differ between replicates, and certain genes do not express in all replications, leading to missing data. Statistical analyses of such problems is an area of increasing concern, and various statistical models and methods have been developed to analyze these situations. Some recent references in this area include [7] and [18], which study the large $p$ small $n$ problem. The reference [15] studies the joint asymptotics in the context of general regression problems when the number of parameters diverge to infinity with the sample size. In particular, [7] investigates the simultaneous estimation of the marginal distributions in the large $p$ small $n$ problem, and it describes how these results can then be used to control the so-called false discovery rates (FDR).

The primary focus of this paper is to develop a general framework for joint statistical analysis of parameters in high-dimensional problems. Furthermore, we allow a random number of parameters and missing data in our data structures. This is achieved using infinite-dimensional techniques. Although the methods of our paper apply generally to many high-dimensional data problems, we will frequently use the terminology from microarrays to facilitate connections to one of the contemporary scientific disciplines. Now we turn to some specifics of our model.

For each fixed integer $n \geq 1$, we begin with a collection of independent sequences of real valued random variables $\{\xi_{n,i,j} : j \geq 1\}$. All are assumed to be defined on a common probability space, and there is no dependence relationship assumed as $n$ and $j$ vary. In the context of microarrays, for $n$ fixed, each of these sequences represents the expression levels of genes in one replication of the experiment. The index $n$ can be interpreted as either the time frame or as a label for the laboratory where the experiment is being performed. In particular, the random variable $\xi_{n,i,j}$ can then be thought of as the expression level of the $j$th gene in the $i$th replicate with index $n$. The number of replicates, for fixed index $n$, could be any integer $r(n)$, but for the sake of simplicity we take $r(n) = n$. Nevertheless, the techniques of this paper can be applied to develop results for other choices of $r(n)$.

Since the expressed genes between replicates may not coincide, either due to the random number that appear or for other reasons (which can be viewed as random deletions), we incorporate these two nonmutually exclusive possibilities into our model. We let $N_{n,i}$ denote the random number of variables within the $i$th replicate having index $n$. We also assume for each integer $n \geq 1$ that $\{N_{n,i} : i \geq 1\}$ is an i.i.d. sequence of integer valued random variables with $P(N_{n,i} \geq 1) = 1$. Of course, in real datasets for fixed $n$, the row lengths are bounded, but our results also apply to situations where they are unbounded.

To model missing data, we postulate that the missing mechanism is independent of the expression level and the random number of parameters involved. For this reason, we introduce the Bernoulli random variables $\{R_{n,i,j} : n \geq 1, i \geq 1, j \geq 1\}$ to represent missing data indicators, where

$$(1.1) \qquad P(R_{n,i,j} = 1) = p \qquad \text{for } n \geq 1, i \geq 1, j \geq 1.$$

We will assume that $0 < p \leq 1$, and also that the sequences $\{\xi_{n,i,j} : n \geq 1, i \geq 1, j \geq 1\}, \{R_{n,i,j} : n \geq 1, i \geq 1, j \geq 1\}$, and $\{N_{n,i} : n \geq 1, i \geq 1\}$ are independent. The case $p = 1$ corresponds to the case that there is no missing data.

In traditional multivariate analysis, such data is typically represented as random vectors in a fixed dimension $d$. However, since we are studying the model in which the dimension of the parameter vector diverges to infinity with the sample size, we represent it as a vector in $R^\infty$, the linear space of all real sequences. That is, we set

$$(1.2) \qquad \mathbf{X}_{n,i} = \sum_{j \geq 1} \xi_{n,i,j} \theta_{n,i,j} \mathbf{e}_j, \qquad i = 1, \ldots, n,$$

where

$$(1.3) \qquad \theta_{n,i,j} = I(j \leq N_{n,i}) R_{n,i,j},$$

for $n \geq 1$, $i \geq 1$, $j \geq 1$, and $\{\mathbf{e}_j : j \geq 1\}$ is the canonical basis for $R^\infty$; that is, $\mathbf{e}_j = \{\delta_{j,k} : k \geq 1\}$ for $j = 1, 2, \ldots$, where $\delta_{j,k} = 1$ for $j = k$ and $0$ if $j \neq k$. In the context of microarrays, the coordinates of the vector $\mathbf{X}_{n,i}$ are thought to be the "normalized expression levels" of genes identified in the $i$th replicate with index $n$. In probabilistic terms, the collection $\mathbf{X}_{n,1}, \mathbf{X}_{n,2}, \ldots, \mathbf{X}_{n,n}$ forms a triangular array of $n$ independent $R^\infty$-valued random vectors. Let $N_n^\star = \max_{1 \leq i \leq n} N_{n,i}$ denote the maximum number of components (columns) in the dataset; or in the context of microarrays, the total number of expressed genes present. If $P(N_n^\star < \infty) = 1$, the components of $\mathbf{X}_{n,i}$, namely $\xi_{n,i,j} \theta_{n,i,j}$, equal $0$ for $j > N_{n,i}$. In other words, $\mathbf{X}_{n,i} \in c_0$, where $c_0$ is the linear space of all real sequences converging to 0. Hence, we will be concerned with asymptotic inference for data in $c_0$. Throughout the paper, we allow the possibility that $P(N_n^\star = \infty) > 0$. We also will use the notation $\mathbf{x} = \sum_{j \geq 1} x_j \mathbf{e}_j$ to denote a typical vector in $R^\infty$, where $\{\mathbf{e}_j : j \geq 1\}$ denotes the canonical basis vectors defined above.

The space $c_0$, with the usual sup-norm given by

$$(1.4) \qquad \|\mathbf{x}\|_\infty = \sup_{i \geq 1} |x_i|,$$

is naturally appropriate when studying the asymptotic inference for a one-sample problem using the maximum of suitable "averages" of gene expressions. In our data analyses, we also use $\ell_\rho$ subspaces, $2 \leq \rho \leq \infty$, determined by the norm

$$\|\mathbf{x}\|_\rho = \left(\sum_{j \geq 1} |x_j|^\rho\right)^{1/\rho}, \tag{1.5}$$

when $2 \leq \rho < \infty$, and by (1.4) when $\rho = \infty$. Related theoretical results for these norms are studied in [8]. Our main asymptotic results concern the statistics

$$\mathbf{S}_{n,n} = \sum_{i=1}^n \mathbf{X}_{n,i} \tag{1.6}$$

and

$$\tilde{\mathbf{S}}_{n,n} = \sum_{i=1}^n \sum_{j \geq 1} \frac{\xi_{n,i,j} \theta_{n,i,j}}{V_{n,j}^{1/2}} \mathbf{e}_j \equiv \sum_{i=1}^n \sum_{j=1}^{N_{n,i}} \frac{\xi_{n,i,j} R_{n,i,j}}{V_{n,j}^{1/2}} \mathbf{e}_j, \tag{1.7}$$

where

$$V_{n,j} = \max\left\{1, \sum_{i=1}^n \theta_{n,i,j}\right\}, \qquad n \geq 1, j \geq 1. \tag{1.8}$$

Here, the coordinate wise random-normalizers $V_{n,j}$ take into account the differences amongst columns due to missing data and random row lengths, and if we replace the $V_{n,j}$ in $\tilde{\mathbf{S}}_{n,n}$ by $n^{1/2}$, then we obtain $\mathbf{S}_{n,n}/n^{1/2}$. Our results include consistency, rates of convergence, and asymptotic normality for these sums. The statistic $\tilde{\mathbf{S}}_{n,n}$ is important when we consider asymptotic normality in our model, as it essentially normalizes each column by the square root of the number of terms in that column. We also study the statistic

$$\tilde{\mathbf{T}}_{n,n} = \sum_{i=1}^n \sum_{j \geq 1} \frac{\xi_{n,i,j} \theta_{n,i,j}}{V_{n,j}} \mathbf{e}_j \equiv \sum_{i=1}^n \sum_{j=1}^{N_{n,i}} \frac{\xi_{n,i,j} R_{n,i,j}}{V_{n,j}} \mathbf{e}_j. \tag{1.9}$$

When $N_{n,i} = p_n$, where $p_n$ is nonrandom, exponential in $n$, and $p = 1$, this is sometimes called the large $p$ small $n$ problem, and [7] and [18] studied the behavior of $\mathbf{S}_{n,n}$ in the sup-norm under various assumptions on the tail behavior of $\xi_{n,i,j}$. For example, the results proved in [18] assume that the random variables have bounded support, while [7] replaces this condition by various exponential decay conditions on the tail behavior of $\xi_{n,i,j}$. The primary technique employed in [7] to obtain consistency results uses the uniform constants for the exponential rate of sup-norm convergence of the empirical distribution function to the true distribution function. This is then used to obtain results for the relevant partial sums of random variables using integration-by-parts techniques. While this approach yields useful results, the integration by parts required seems to obscure the true nature of the matter. From what we do here, we will see that it is more fruitful to study the problem from the point of view of the random variables themselves in that we are able to clarify some of the results described in [7], and also extend them under a broader range of conditions to our more general model. While we focus on the mean functional, [7] studies other interesting functionals of the data.

The rest of the paper is organized as follows. Section 2 presents the main results. These results concern joint consistency and joint asymptotic normality. Section 3 contains applications to hypothesis tests and Section 4 is devoted to simulation results and real data analysis. Section 5 contains the necessary probability estimates, while Sections 6 and 7 contain the proofs of our main results.

Throughout the paper, $Lx = L(x) = \log_e(\max(x, e))$.

**2. The main results.** The results that we obtain will depend critically on the tail probabilities of the random variable $\{\xi_{n,i,j}\}$. These assumptions are of two types, namely that the tail probabilities decay at an exponential rate, or that they decay polynomially. In the large $p$ small $n$ problem, our results imply that these tail probability conditions are closely tied to the way $p$ must relate to $n$. For example, in the classic version of this problem where $p$ grows exponentially fast in $n$, we need tail probabilities that decay exponentially fast, whereas if $p$ grows only as a power of $n$, then we only need polynomial decay for the tails. The precise nature of this interplay for consistency results is contained in Theorems 2.2 and 2.3. In particular, the remarks following these theorems contain precise information on their relationship to the large $p$ small $n$ problem.

First, we discuss the exponential decay case. Here, we assume that for some $r$, $0 < r \leq 2$, and all $x \geq 0$ there are constants $c_{n,j}$ and $k_{n,j}$ such that

$$(2.1) \qquad P(|\xi_{n,i,j}| \geq x) \leq c_{n,j} e^{-k_{n,j} x^r}$$

for all $n \geq 1, j \geq 1$. Random variables satisfying (2.1) with $r = 2$ are usually said to be sub-Gaussian, and if for $1 \leq i \leq n$ we have that each $\xi_{n,i,j}$ takes values in the interval $[a_{n,j}, b_{n,j}]$, then we will see below that (2.1) holds with $r = 2$,

$$(2.2) \qquad c_{n,j} = 2 \quad \text{and} \quad k_{n,j} = (2(b_{n,j} - a_{n,j})^2)^{-1}, \qquad n \geq 1, j \geq 1.$$

Throughout, when $\xi_{n,i,j} = 0$ with probability one, in (2.1) we take

$$(2.3) \qquad c_{n,j} = 1 \quad \text{and} \quad k_{n,j} = \infty, \qquad n \geq 1, j \geq 1.$$

In addition, note that $c_{n,j} \geq 1$ is necessary by setting $x = 0$ in (2.1).

It is also useful to notice that if the condition (2.1) holds for some $r^* > 1$, then it holds for all $1 \leq r \leq r^*$ by simply adjusting the constants $c_{n,j}$ and keeping the same $k_{n,j}$. In particular, if (2.1) holds for some $r > 2$, then it holds for $r = 2$, and we are in the sub-Gaussian setting. In [7], this seems to have gone unnoticed, and there one finds results for $r > 2$ which are weaker than the corresponding $r = 2$ results. However, this should not be the case as the previous comment implies the $r = 2$ result applies directly to what is proved there. Of course, in some settings there could be results that distinguish between various $r$ values, even for $r > 2$, but that does not happen here, and is why we restrict $r$ to be in $(0, 2]$. Our methods also yield results when $0 < r < 1$, whereas in [7], the parameter $r$ is always greater than equal to one.

Another situation we will discuss is when the assumption of exponential decay of the tails of $\xi_{n,i,j}$ in (2.1) is replaced by the polynomial decay

$$(2.4) \qquad P(|\xi_{n,i,j}| \geq x) \leq \frac{c_{n,j}}{(1+x)^{k_{n,j}}}, \qquad x \geq 0,$$

where $c_{n,j} \geq 1$ and typically for our results, $2 < k_{n,j} < \infty$.

We will assume throughout the paper that $E(\xi_{n,i,j}) = 0$ for all $n, i, j \geq 1$. Should this not be the case, one would simply replace the tail probability conditions in (2.1) and (2.4) by analogous conditions for the variables $\{\xi_{n,i,j} - E(\xi_{n,i,j})\}$, and formulate the results in terms of these variables.

2.1. *Consistency and rates of convergence.* In this subsection we present several consistency and rate of convergence results for $\mathbf{S}_{n,n}$ and $\tilde{\mathbf{S}}_{n,n}$.

THEOREM 2.1. *Let $\{\mathbf{X}_{n,i} : 1 \leq i \leq n\}$ be as in (1.2), assume (2.1) holds with $r = 2$, and take $\{a_n : n \geq 1\}$ to be a sequence of positive numbers. Furthermore, assume $c_{n,j}, k_{n,j}$ are constants such that $c_{n,j} \geq 1, k_{n,j} \leq \infty$ and*

$$(2.5) \qquad \sum_{n \geq 1} \sum_{j \geq 1} \exp\{-(\varepsilon a_n)^2 k_{n,j}/(16 c_{n,j})\} < \infty$$

*for all $\varepsilon > \varepsilon_0$. Then*

$$(2.6) \qquad \sum_{n \geq 1} P(\|\tilde{\mathbf{S}}_{n,n}\|_\infty \geq \varepsilon a_n) < \infty$$

*for all $\varepsilon > \varepsilon_0$, where $\tilde{\mathbf{S}}_{n,n}$ is given as in (1.7). Thus, if the constants $c_{n,j}$ and $k_{n,j}$ are such that uniformly in $n \geq 1$ and for some $\delta > 0$,*

$$(2.7) \qquad k_{n,j}/(16 c_{n,j}) \geq \delta L(j+3),$$

*then for all $\varepsilon > 0$ such that $\varepsilon^2 \delta > 1$ we have*

$$\limsup_{n \to \infty} \frac{\|\tilde{\mathbf{S}}_{n,n}\|_\infty}{(L(n+3))^{1/2}} \leq \varepsilon. \tag{2.8}$$

*In particular, if (2.7) holds, then with probability one*

$$M = \sup_{n \geq 1} \frac{\|\tilde{\mathbf{S}}_{n,n}\|_\infty}{(L(n+3))^{1/2}} < \infty \tag{2.9}$$

*and there exists an $\alpha > 0$ such that*

$$E(e^{\alpha M^2}) < \infty. \tag{2.10}$$

*Moreover, if the $V_{n,j}^{1/2}$ are replaced by $n^{1/2}$ in $\tilde{\mathbf{S}}_{n,n}$, then again (2.6), (2.8) and (2.10) continue to hold.*

REMARK 2.1. In Theorem 2.4, we will establish a central limit theorem for $\tilde{\mathbf{S}}_{n,n}$, and that $\tilde{\mathbf{T}}_{n,n}$ converges to zero in probability under related conditions.

In Theorem 2.1, the impact of the random row sizes $\{N_{n,i} : i \geq 1\}$ is hidden due to our choice of normalizations $\{a_n\}$ as given in (2.5). For example, (2.5) implies the ratio $k_{n,j}/c_{n,j}$ cannot be bounded as $j$ goes to infinity, but in our next result we only require this ratio to be uniformly bounded below in both $n$ and $j$ by a strictly positive constant. Under this different set of conditions, the role of $\{N_{n,i} : i \geq 1\}$ appears in the normalizations for $\tilde{\mathbf{S}}_{n,n}$ given by $h(n)$ in (2.12). In particular, if $N_{n,i} = p_n$ for $\{i \geq 1, n \geq 1\}$, then Theorem 2.2 and Remark 2.4 below yield the results in Corollaries 1 and 2 in [7] when $r = 2$. The consistency results in [7] for $1 \leq r < 2$, as well as for many other cases, follow immediately from Theorem 2.3 below.

THEOREM 2.2. *Let $\{\mathbf{X}_{n,i} : 1 \leq i \leq n\}$ be as in (1.2), and assume (2.1) holds with $r = 2$, and that for $1 \leq c < \infty, 0 < k < \infty$, we have $c_{n,j} \leq c$ and $k_{n,j} \geq k$ for all $n, j \geq 1$. Let*

$$h(n) = (\theta_1^{-1} L(E(N_n^*)) + \theta_2 L(n))^{1/2}, \tag{2.11}$$

*where $\theta_1 = k/(16c)$ and $\theta_2 > 0$. Then*

$$\lim_{n \to \infty} P\left( \frac{\|\tilde{\mathbf{S}}_{n,n}\|_\infty}{h(n)} \geq 1 \right) = 0 \tag{2.12}$$

*and if also $\theta_1 \theta_2 > 1$, then*

$$\sum_{n \geq 1} P(\|\tilde{\mathbf{S}}_{n,n}\|_\infty \geq h(n)) < \infty. \tag{2.13}$$

*Finally, if the $V_{n,j}^{1/2}$ are replaced by $n^{1/2}$ in $\tilde{\mathbf{S}}_{n,n}$, then again (2.12) and (2.13) hold.*

REMARK 2.2. Note that (2.12) immediately implies

$$\frac{\|\tilde{\mathbf{S}}_{n,n}\|_\infty}{n^{1/2}} = O_P\left(\left(\frac{L(E(N_n^*)) + L(n)}{n}\right)^{1/2}\right). \quad (2.14)$$

In particular, if $L(E(N_n^*))/n$ converges to zero, then (2.14) implies that $\|\tilde{\mathbf{S}}_{n,n}\|_\infty/n^{1/2}$ tends to zero in probability. In addition, (2.13) implies with probability one that

$$\limsup_{n\to\infty} \frac{\|\tilde{\mathbf{S}}_{n,n}\|_\infty}{h(n)} \leq 1 \quad (2.15)$$

and hence $\tilde{\mathbf{S}}_{n,n}/n^{1/2}$ converges to zero with probability one provided $\theta_1\theta_2 > 1$ and $\lim_{n\to\infty} n^{-1}L(E(N_n^*)) = 0$. Furthermore, if $N_{n,i} = p_n \geq n$, $\{i \geq 1,\ n \geq 1\}$, then (2.14) immediately relates to the results of Corollary 2 in [7], as it implies

$$\frac{\|\tilde{\mathbf{S}}_{n,n}\|_\infty}{n^{1/2}} = O_P\left(\left(\frac{L(p_n)}{n}\right)^{1/2}\right). \quad (2.16)$$

In particular, (2.14) improves Corollary 2 and its proof considerably whenever $r \geq 2$ there, and the case $0 < r < 2$ will be discussed in what follows. Once we establish Lemma 1 below these results also apply to Corollary 1 of [7] in a standard way.

We next study the situation when the random variables $\{\xi_{n,i,j}\}$ satisfy the exponential tail condition (2.1) with $0 < r < 2$, or polynomial decay as in (2.4). When $1 \leq r < 2$, a special case of these results clarifies Corollary 2 of [7]. This can be seen in Remark 2.3 below. The $r = 2$ case in this corollary already appeared in (2.16) when $V_{n,j} = n$. It should also be observed that Theorem 2.3 provides sufficient conditions for consistency which involve a precise relationship between the size of $p_n$ in the large $p$ small $n$ problem, and the tail decay of the data. This relationship is shown to exist even when there is only polynomial decay in the data, and as one might expect in this situation the growth of $p_n$, or $E(N_n^*)$, needs to be further restricted, that is, in such results $p_n$ and $E(N_n^*)$ grow at a corresponding polynomial rate.

THEOREM 2.3. *Let $\{\mathbf{X}_{n,i} : 1 \leq i \leq n\}$ be as in (1.2) and assume that (2.1) holds with $0 < r < 2$. Also assume for all $n \geq 1$ and $j \geq 1$, that $c_{n,j} \leq c$ and $k_{n,j} \geq k$, where $1 \leq c < \infty, 0 < k < \infty$. Let $s_n = c_1(L(E(N_n^*)) + 2L(n))^{1/r}$, and*

$$h(n) = (c_2^{-1} L(E(N_n^*)) + c_3 L(n))^{1/2}, \quad (2.17)$$

*where $c_1 > 2/k^{1/r}$, $c_2 = k/(128c)$, and $c_3 > 0$. Then*

$$\lim_{n\to\infty} P\left(\frac{\|\mathbf{S}_{n,n}\|_\infty}{n^{1/2}s_n h(n)} \geq 1\right) = 0. \quad (2.18)$$

*If we also assume $c_2 c_3 > 1$, then*

$$(2.19) \qquad \sum_{n \geq 1} P(\|\mathbf{S}_{n,n}\|_\infty \geq n^{1/2} s_n h(n)) < \infty.$$

*Furthermore, if $k > 2$ and the polynomial condition in (2.4) holds, then*

$$(2.20) \qquad \|\mathbf{S}_{n,n}\|_\infty = O_P(n^{1/2} s_n (L(E(N_n^*)))^{1/2}),$$

*where $s_n = (nE(N_n^*))^{1/k+\beta}$ and $\beta > 0$. Additionally, if $E(N_n^*) \geq n$, $b > 8$, and $k\beta > 1/2$, then*

$$(2.21) \qquad \sum_{n \geq 1} P(\|\mathbf{S}_{n,n}\|_\infty \geq b s_n n^{1/2} (L(E(N_n^*)))^{1/2}) < \infty.$$

*In particular, if $E(N_n^*)$ is asymptotic to $n^\gamma$ for $\gamma \geq 1$, then*

$$(2.22) \qquad \sum_{n \geq 1} P(\|\mathbf{S}_{n,n}\|_\infty \geq b s_n n^{1/2} (L(E(N_n^*)))^{1/2}) < \infty,$$

*provided $b > 8$ and $(\gamma + 1)k\beta > 1$.*

REMARK 2.3. An immediate consequence of (2.18) is that

$$(2.23) \qquad \frac{\|\mathbf{S}_{n,n}\|_\infty}{n} = O_P\left(\frac{(L(E(N_n^*)) + L(n))^{(2+r)/(2r)}}{n^{1/2}}\right)$$

and if $(L(E(N_n^\star)))^{(2+r)/(2r)}/n^{1/2} \to 0$, then (2.23) easily implies $\mathbf{S}_{n,n}/n$ converges to zero in probability. In addition, if $N_{n,i} = p_n$ for $n \geq 1, i \geq 1$, where $\{p_n : n \geq 1\}$ is a sequence of integers, and $p_n \geq n$, then it follows from (2.23) that

$$(2.24) \qquad \frac{\|\mathbf{S}_{n,n}\|_\infty}{n} = O_P\left(\frac{(L(p_n))^{(2+r)/(2r)}}{n^{1/2}}\right).$$

Hence, using the above for $r \in (0, 2)$, and (2.16) for the case $r = 2$, one obtains an extension and clarification of Corollary 2 and its proof in [7]. It is also interesting to observe that the method of proof for Theorem 2.3 applied to the $r = 2$ situation only yields

$$(2.25) \qquad \frac{\|\mathbf{S}_{n,n}\|_\infty}{n} = O_P\left(\frac{L(p_n)}{n^{1/2}}\right).$$

Hence, we see the methods used for the $r = 2$ case in Theorem 2.2 are sharper than those we have for other values of $r$.

REMARK 2.4. Under the assumption of polynomial decay given in Theorem 2.3, and assuming that $E(N_n^*)$ is asymptotic to $n^\gamma$ for $\gamma \geq 1$, we easily see from (2.22) that $\|\mathbf{S}_{n,n}\|_\infty/n$ converges to zero almost surely provided $k$ is sufficiently large so that for $\beta > 0$ we have $(\gamma + 1)\beta k > 1$ and $(\gamma + 1)(1/k + \beta) < 1/2$.

2.2. *Asymptotic normality results.* In this section, we present results on the asymptotic normality of the quantity $\tilde{\mathbf{S}}_{n,n}$. Since this estimator typically lives in $c_0$, Theorem 2.4 is a central limit theorem in that setting. Nevertheless, we also have proved CLTs in $\ell_\rho, 2 \leq \rho < \infty$, similar to that found in Theorem 2.4. They appear in [8]. These results hold when the underlying process is a triangular array with random row lengths and possibly missing data. We also are able to use the coordinate-wise random normalizations $V_{n,j}^{1/2}$. However, we also use classical normalizations for some of the $\ell_\rho$ results, and in that case the related CLTs hold under far weaker moment conditions. See [1], page 206, for some classical results. The paper [14] contains CLTs in $c_0$, as well as related references, and much is known about the CLT in the spaces $\ell_\rho, 2 \leq \rho < \infty$. However, none of these results incorporate random row lengths, missing data, or coordinate-wise random normalizations in their formulations. In addition, the results in [14] require a uniform boundedness assumption on the $\{\xi_{n,i,j}\}$ to obtain results related to what we prove. Finally, we mention that in our simulation results we include the use of our CLTs in $\ell_\rho$, $2 \leq \rho < \infty$.

A key assumption in any central limit theorem is that there is a limiting covariance function. Since our results include the use of random column-wise normalizers, we have need of a couple different limiting covariances. That is, if

$$\Gamma_n(j_1, j_2) = \sum_{i=1}^{n} E(\xi_{n,i,j_1} \xi_{n,i,j_2})/n$$

is such that

(2.26) $$\lim_{n \to \infty} \Gamma_n(j_1, j_2) = \Gamma(j_1, j_2)$$

for all $j_1, j_2 \geq 1$, then for $k = 1, 2$ we set

(2.27) $$\Gamma(k, j_1, j_2) = p^k \Gamma(j_1, j_2) \qquad \text{for } j_1 \neq j_2,$$

and

(2.28) $$\Gamma(k, j_1, j_2) = p^{k-1} \Gamma(j_1, j_2) \qquad \text{for } j_1 = j_2.$$

THEOREM 2.4. *Let $\{\mathbf{X}_{n,i} : 1 \leq i \leq n\}$ be as in (1.2), assume (2.1) holds with $r = 2$, and that $c_{n,j}, k_{n,j}$ are constants such that $c_{n,j} \geq 1, k_{n,j} < \infty$ and*

(2.29) $$\sup_{n,j \geq 1} c_{n,j}/k_{n,j} < \infty.$$

*Also assume for all $\delta > 0$ that*

(2.30) $$\lim_{d \to \infty} \sup_{n \geq 1} \sum_{j \geq d} \exp\{-\delta k_{n,j}/c_{n,j}\} = 0.$$

If $\tilde{\mathbf{S}}_{n,n}$ is given as in (1.7), then

(2.31)  $\quad\quad\quad\quad\quad \{\mathcal{L}(\tilde{\mathbf{S}}_{n,n}) : n \geq 1\}$ is tight in $c_0$.

In addition, if $\tilde{\mathbf{T}}_{n,n}$ is as in (1.9), and for each $j \geq 1$ we have $\lim_{n\to\infty} P(N_{n,1} < j) = 0$, then $\tilde{\mathbf{T}}_{n,n}$ converges in probability to zero in $c_0$. Moreover, if the $V_{n,j}^{1/2}$ are replaced by $n^{1/2}$ in $\tilde{\mathbf{S}}_{n,n}$, then again (2.31) holds, and $\tilde{\mathbf{T}}_{n,n}$ converges in probability to zero. Furthermore, if we also assume (2.26), (2.27), and (2.28) hold, and for each $d < \infty$ we have $P(\min_{1 \leq i \leq n} N_{n,i} < d) = o(1/n^2)$ as $n$ tends to infinity, then $\Gamma(k, \cdot, \cdot)$ is the covariance of a centered Gaussian measure $\gamma_k$ on $c_0$ for $k = 1, 2$, and

(2.32)  $\quad\quad\quad\quad\quad \mathcal{L}(\tilde{\mathbf{S}}_{n,n})$ converges weakly to $\gamma_1$

on $c_0$. If the $V_{n,j}^{1/2}$ are replaced by $n^{1/2}$, then (2.32) still holds with limiting measure $\gamma_2$.

REMARK 2.5. The conditions (2.26), (2.29), and (2.30), along with (2.1) when $r = 2$, allow the limiting Gaussian measures $\gamma_k$ to exist on $c_0$. Moreover, without such assumptions, with the most important being (2.26) and (2.30), there are examples of triangular arrays of the form indicated when the CLT must fail on $c_0$, although it may hold on $R^\infty$. Of course, without (2.26), then the CLT will fail even on $R^\infty$.

**3. Application to hypothesis tests.** In this section, we deal with application of our results to one-sample problem, two-sample problem, and one-way random effects models. Joint hypothesis testing was also considered in [11]. We assume throughout this section that the distribution of the random variables $\mathbf{X}_{n,i}$ in (1.2) are independent of $i$. This implies that $E(\xi_{n,i,j}) = \mu_{n,i,j}$ is independent of $i$ and we write it as $\mu_{n,j}$.

3.1. *One-sample and two-sample problems.* In this subsection, we apply our results to test if the "mean vector" equals a specified vector in the one-sample case and if the difference in the "mean vectors" is zero in the two-sample case. More precisely, for the one-sample case consider testing the null hypothesis $\mathbf{H}_0 : \boldsymbol{\mu}_n = \mathbf{0}$, where $\boldsymbol{\mu}_n$ is an infinite-dimensional vector whose components are $\mu_{n,j}$. The quantity $\tilde{\mathbf{S}}_{n,n}$, defined by (1.7), can be used for developing a test of $\mathbf{H}_0$. To this end, let us denote the data vectors by $\vec{\mathbf{X}}_n = \{\mathbf{X}_{n,1}, \ldots, \mathbf{X}_{n,n}\}$. One can use the $\ell_\rho$ norm for $\rho \geq 2$ and the $c_0$ norm to define various nonrandomized test functions $\phi_\rho(\vec{\mathbf{X}}_n)$ as follows:

(3.1)  $\quad\quad\quad\quad\quad \phi_\rho(\vec{\mathbf{X}}_n) = \begin{cases} 1, & \text{if } \|\tilde{\mathbf{S}}_{n,n}\|_\rho > c, \\ 0, & \text{otherwise,} \end{cases}$

where $c = c_\rho$ is so chosen that $E(\phi_\rho(\vec{\mathbf{X}}_n)|\mathbf{H}_0) \leq \alpha$. The test function based on the $c_0$ norm is given by

$$\phi_\infty(\vec{\mathbf{X}}_n) = \begin{cases} 1, & \text{if } \|\tilde{\mathbf{S}}_{n,n}\|_\infty > c, \\ 0, & \text{otherwise,} \end{cases} \tag{3.2}$$

where $c = c_\infty$ is so chosen that $E(\phi_\infty(\vec{\mathbf{X}}_n)|\mathbf{H}_0) \leq \alpha$. In the context of the two-sample problem, the null hypothesis is $\mathbf{H}_0 : \boldsymbol{\mu}_n^1 = \boldsymbol{\mu}_n^2$, where $\boldsymbol{\mu}_n^k$ represents the infinite-dimensional mean vector from the $k$th population. Now, using $\tilde{\mathbf{S}}_{n,n}^{(k)}$ to denote $\tilde{\mathbf{S}}_{n,n}$ for the $k$th population [note that these are constructed using a superscript $k$ in all the quantities in (1.2) and (1.3) and these quantities are independent in $k$], the test function for the $\ell_\rho$ norm is

$$\phi_\rho^{(2)}(\vec{\mathbf{X}}_n) = \begin{cases} 1, & \text{if } \|\tilde{\mathbf{S}}_{n,n}^{(1)} - \tilde{\mathbf{S}}_{n,n}^{(2)}\|_\rho > c, \\ 0, & \text{otherwise,} \end{cases} \tag{3.3}$$

where $c = c_\rho$ is so chosen that $E(\phi_\rho^{(2)}(\vec{\mathbf{X}}_n)|\mathbf{H}_0) \leq \alpha$. The test function based on the $c_0$ norm is given by

$$\phi_\infty^{(2)}(\vec{\mathbf{X}}_n) = \begin{cases} 1, & \text{if } \|\tilde{\mathbf{S}}_{n,n}^{(1)} - \tilde{\mathbf{S}}_{n,n}^{(2)}\|_\infty > c, \\ 0, & \text{otherwise,} \end{cases} \tag{3.4}$$

where $c = c_\infty$ is so chosen that $E(\phi_\infty^{(2)}(\vec{\mathbf{X}}_n)|\mathbf{H}_0) \leq \alpha$. We observe that under our formulation unequal sample sizes from the two populations are allowed. In some applications, the number of components in the two groups and the sample sizes coincide. In these cases, as in traditional multivariate analysis, the two-sample problem reduces to the one-sample problem.

To perform the test one needs the distribution of $\|\tilde{\mathbf{S}}_{n,n}\|_\rho$, $\|\tilde{\mathbf{S}}_{n,n}\|_\infty$, $\|\tilde{\mathbf{S}}_{n,n}^{(1)} - \tilde{\mathbf{S}}_{n,n}^{(2)}\|_\rho$, and $\|\tilde{\mathbf{S}}_{n,n}^{(1)} - \tilde{\mathbf{S}}_{n,n}^{(2)}\|_\infty$. The following proposition provides the asymptotic distribution for these statistics. It's proof can be obtained as a corollary of Theorem 2.4 for the $c_0$ case, and Theorem 7 and Remark 11 of [8] for $\ell_\rho$ spaces, and using the fact that the distribution function of the norm with respect to a Gaussian measure on a separable Banach space is continuous.

PROPOSITION 3.1. *If the appropriate null hypothesis holds, then Theorem 2.4 implies that $P(\|\tilde{\mathbf{S}}_{n,n}\|_\infty > c)$ and $P(\|\tilde{\mathbf{S}}_{n,n}^{(1)} - \tilde{\mathbf{S}}_{n,n}^{(2)}\|_\infty > c)$ converges to $P(\|\mathbf{G}\|_\infty > c)$ and $P(\|\mathbf{G}^{(1)} - \mathbf{G}^{(2)}\|_\infty > c)$, respectively, for all $c > 0$, where $\mathcal{L}(\mathbf{G}) = \mathcal{L}(\mathbf{G}^{(1)}) = \mathcal{L}(\mathbf{G}^{(2)}) = \gamma$ and $\gamma$ is the Gaussian measure identified there. Furthermore, $\mathbf{G}^{(1)}$ and $\mathbf{G}^{(2)}$ are independent. A similar result holds under the conditions of Theorem 7 of [8] when the infinity norm is replaced by the $\rho$-norm.*

3.2. *One-way random effects models.* In the analysis of gene-expression data, the random effects model with random plate effect is often used for data analysis. That is, if $\xi_{n,i,j,k}$ are the expression levels of the $j$th gene in the $i$th replicate, receiving the $k$th treatment, then for $k = 1, 2$ the one-way random effects model is given by

$$\xi_{n,i,j,k} = \mu_{n,j,k} + T_{n,i,k} + \varepsilon_{n,i,j,k},$$
(3.5)
$$i = 1, 2, \ldots, n, j = 1, 2, \ldots, b(n),$$

where $\mu_{n,j,k}$ is the mean expression level for the $j$th gene receiving the $k$th treatment in the lab $n$, $T_{n,i,k}$ are independent Gaussian random variables with 0 and variance $\sigma_{n,k}^2$ and $\varepsilon_{n,i,j,k}$ are i.i.d. random variables with mean 0 and variance $\sigma^2$. More complicated models that take into account tip effect and dye effects have been studied in the applied literature. The index $n$ in the subscript is usually suppressed in the applied literature but we keep it to show the relationship with our model. The random variables $T_{n,i,k}$ introduce correlations in the expression levels of $\xi_{n,i,j,k}$ across genes and a standard calculation shows that this correlation structure is compound symmetric. This model can be seen as a particular case of our model with $N_{n,i} = b(n)$, $R_{n,i,j} \equiv 1$, and compound symmetric covariance matrix. Proposition 3.1 above can be used for performing hypothesis tests concerning the expression levels of multiple genes simultaneously. Furthermore, the models developed in the paper allow for some extensions of the one-way random effects models to incorporate missing data and random number of parameters.

**4. Simulation results and real data analysis.**

4.1. *Simulation results.* In this section, we evaluate our methodology, using simulations, when the number of replications is small, but the number of variables is large. All our simulation results are based on 5000 independent trials of 10 replications. We purposely chose $n$ small to reflect many real applications.

As a first step, we need to "approximate" the limiting distribution of the random variables appearing in Proposition 3.1. We will work with the case $N_{n,i} = b(n)$, and assume that $\mathbf{X}_{n,1}, \ldots, \mathbf{X}_{n,n}$ are $n$ i.i.d. $b(n)$-dimensional vectors with distribution $G_n(\cdot)$, whose coordinates have tails that satisfy the sub-Gaussian property. Let $\hat{\Sigma}_n$ denote an estimate of the covariance matrix $\Sigma_n = ((\sigma_{n,u,v}))$, where $\Sigma_n$ is a $b(n) \times b(n)$ matrix given by

(4.1) $$\sigma_{n,u,v} = E(\xi_{n,1,u} - \mu_{n,u}^{(0)})(\xi_{n,1,v} - \mu_{n,v}^{(0)}).$$

In the above definition, $\mu_{n,u}^{(0)}$ and $\mu_{n,v}^{(0)}$ are the specified values under the null hypothesis. Note that $\hat{\Sigma}_n$ is a function of the data vector $\vec{\mathbf{X}}_{n,n}$. One choice for $\hat{\Sigma}_n$ is the sample covariance matrix. In fact, better options are available, and we will explain them later below. If $\hat{\Sigma}_n$ is positive definite, then given $\vec{\mathbf{X}}_{n,n}$, we generate $t$ i.i.d. random vectors $\mathbf{Y}_{n,i}$ of dimension $b(n)$ whose distribution is Gaussian with mean vector $\mathbf{0}$ and covariance matrix $\hat{\Sigma}_n$; that is,

$$(4.2) \qquad \mathbf{Y}_{n,i}|\vec{\mathbf{X}}_n \sim N_{b(n)}(\mathbf{0}, \hat{\Sigma}_n) \qquad \text{a.s.}, 1 \leq i \leq t.$$

We will call $\mathbf{Y}_{n,i}$ the Monte Carlo (MC) samples, and throughout the simulations $t = 2000$. We will use $\|\cdot\|$ to denote the $l_\rho$ norm ($\rho \geq 2$) or the $c_0$ norm depending on the space being used. Let $\|\mathbf{Y}_{n,1}\|, \ldots, \|\mathbf{Y}_{n,n}\|$ denote the norms of the MC samples. Furthermore, consider the following nonparametric density estimator; namely, for $x \in R$,

$$(4.3) \qquad h_t(x) = \frac{1}{tc_t} \sum_{i=1}^{t} K\left(\frac{x - \|\mathbf{Y}_{n,i}\|}{c_t}\right),$$

where $c_t$ is a sequence of positive constants converging to 0 such that $tc_t \to \infty$, and $K(\cdot)$ is a density function with $\int_R tK(t)\,dt = 0$. In the above, we have suppressed the dependence on $n$ and on $\omega$ since $n$ and $\omega$ will be held fixed in this discussion. It follows from Devroye [2] that as $t \to \infty$, that for every fixed $\omega \in \Omega$, $h_t(x)$ converges almost everywhere with respect to Lebesgue measure and in $L_1$ to the probability density of the random variable $\|N(\mathbf{0}, \hat{\Sigma}_n)\|$. In all our numerical experiments, we will take $K(\cdot)$ to be a standard normal density, $t = 2000$, and fix the window width $c_t$ at 0.7. Figure 1(a) presents the graph of the density function for the 2-norm, the 10-norm, and the sup-norm. We use these densities to "approximate" the tail probabilities of the norms of the limiting Gaussian appearing in Proposition 3.1. Figure 1(b) shows that the $p$-values from these hypothesis tests are uniformly distributed when we use the sup-norm. This histogram was generated using the $p$-values from a hypothesis test for data generated using a compound symmetric covariance structure. The vertical axis represents the proportion of times the $p$-values belonged to a particular range.

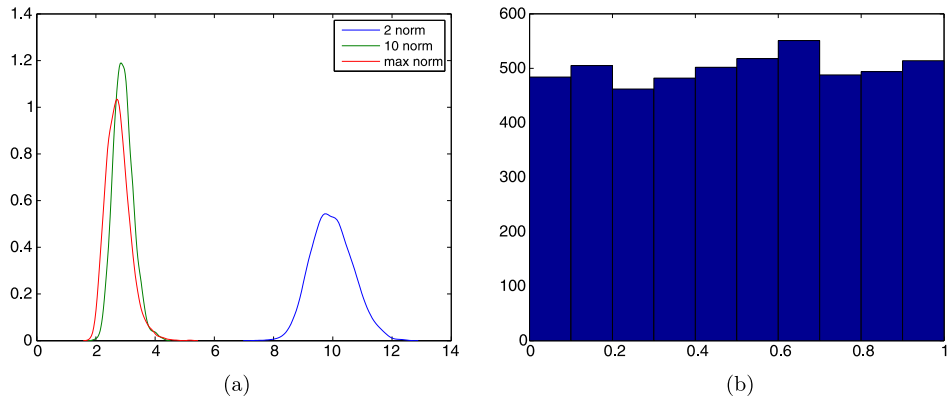

Fig. 1. *Kernel density estimates of the norms of statistics and the histogram of p-values under the sup norm.*

*Structured and unstructured covariance estimation, shrinkage and sparsity.* The methodology described above requires an estimate of the covariance matrix. It is folklore in statistics literature that covariance matrix estimation is a hard problem, and the difficulties increase when the number of variables is much larger than the sample size. The papers [10] and [17], as well as others, have clearly demonstrated that the sample covariance matrix behaves poorly in terms of the mean square error. This difficulty is frequently due to having a large number of parameters to estimate, and a limited number of observations available to estimate them. Hence, it is reasonable to expect that if one uses structured covariances some of these difficulties can be mitigated. However, one has to be cautious since it is also well known that assuming independence when correlations are present lead to substantial bias in type-I error rates. Extensive simulations concerning these issues are described in detail in [8] and [9]. Thus, to make the methodology presented here useful and applicable, we need to describe how to handle the covariance matrix estimation.

In studies involving the joint analysis of multiple cDNA microarray data, [16], one frequently encounters covariances that have a block structure. This phenomenon also occurs in the context of sparse covariance estimation where regularization is adopted. Borrowing the idea from shrinkage estimation, [10] developed an estimator of $\Sigma_n$ by taking a convex combination of the unstructured sample covariance matrix and a structured covariance matrix. Their estimator is given by

$$(4.4) \qquad \Sigma_n^\star = (1-\lambda)\hat{\Sigma}_n + \lambda \tilde{\Sigma}_n,$$

where $\hat{\Sigma}_n$ is the method of moments estimator and $\tilde{\Sigma}_n$ is an estimator assuming a particular structure for the covariance matrix. The parameter $\lambda$ can be estimated from the data, and has a closed form expression when $\tilde{\Sigma}_n$ is taken to be the identity matrix, or is compound symmetric, heterogeneous compound symmetric, as well as many other structures. We propose to use (4.4) as the estimator of the covariance matrix needed for generating samples in (4.2). Our extensive simulations described in [8] show that when one shrinks the variances and the covariances, the type-I error rate approaches the nominal values for the test based on the sup-norm. These methods can also be applied to situations when there is mean sparsity. In these cases, one can apply the LASSO type algorithm for estimation purposes and then use our methodology for hypothesis testing.

TABLE 1
*Type I error rates with unstructured covariance matrix and shrinkage*

|              | $\rho = 2$ | $\rho = 4$ | $\rho = \infty$ |
|---|---|---|---|
| $b(n) = 100$  | 0.1312 | 0.0838 | 0.0686 |
| $b(n) = 500$  | 0.1914 | 0.0874 | 0.057  |
| $b(n) = 1000$ | 0.237  | 0.1042 | 0.0526 |

*Simulation analysis using real data information.* We now describe a numerical experiment which examines whether some of the difficulties described in the previous subsection due to covariance matrix estimation can be minimized by using the shrinkage method for estimating the covariance matrix. To make the model closely resemble micro-array data, we considered the leukemia data set described in [6] which studies the gene expression in two types of leukemia, acute lymphoblastic leukemia (ALL) and acute myeloid leukemia (AML). We use the same preprocessing step as described in [3], Section 3.1; retaining 3571 genes from 72 patients, 38 ALL and 25 AML. We apply the standardization technique described in Section 3.3 of [3] on these retained genes. For our simulations, we consider information from the AML group only. The nominal type-I error rate is taken to be 5% in all of our experiments.

The simulation experiments are based on data generated from a $b(n)$-dimensional normal random variable with $n = 10$. We now describe how the mean and covariance matrix are obtained from the data. First, we fix a $b(n)$ and randomly generate $b(n)$ genes from the 3571 genes. From the 25 AML patients, we estimate the $b(n)$-dimensional mean vector $\boldsymbol{\mu}_n$ and the corresponding covariance matrix $\Sigma_n$ by shrinking the covariances as described in [17]. The shrinkage method yields a nonsingular covariance matrix. For all the 5000 simulations, this mean vector and the covariance matrix are fixed. Now, using the shrinkage method described above, we apply our methodology. The resulting type-I error rate is described in the Table 1 above.

From the table, we notice that the type-I error rate is closer to the nominal value for larger values of $\rho$. Further analysis (data not presented here) shows that this phenomenon existed for the other choices of variances and covariances as well (see [8]).

4.2. *Real data analysis.* In a typical microarray situation, one is interested in identifying if a set of genes are differentially expressed. This is a two-sample problem and one can use the methods described in Proposition 3.1 to address this problem. We analyzed the leukemia data set described in [6], which also was used in the previous section. We analyzed the sixteen genes given in Table 2 of [19]. The value of the test statistics defined in (3.3) were determined to be 4.4992, 2.7396 for $\rho = 2$ and $\rho = 4$, respectively. Also the value of the test statistic in (3.4) was determined to be 2.0250. The covariance matrices needed to apply the methodology were calculated by shrinking the variances and covariances using the algorithm described in [17]. The $p$-values corresponding to these test statistics were all less than $10^{-6}$ showing that the genes are differentially expressed. The same conclusion was also obtained by Yan et al. [19] using different methods. More importantly, as explained in [19], these genes have biological significance and the three existing statistical methods in popular use did not identify them to be differentially expressed. Even though this example involved only sixteen genes, the simulation results in the previous section show that such a limitation is not necessary. Thus, this analysis combined with our simulation results described above, show the importance and usefulness of the proposed methodology.

**5. Some probability estimates.** Here, we provide some basic probability estimates used throughout the paper. The first lemma deals with the sub-Gaussian situation, and the inequality we present for bounded random variables is not best possible, as slightly better constants in the basic estimate can be obtained from Theorem 1 in [5]. Nevertheless, we include a proof for this case, as our argument generalizes to the unbounded case. Our approach is to compute the necessary Laplace transforms, and then use Markov's inequality efficiently. This is standard for such problems, but in order to proceed from first principles, and also keep track of relevant constants, we include the details.

LEMMA 5.1. *Let $X_1, \ldots, X_n$ be independent random variables with $E(X_i) = 0$. If $P(X_i \in [a, b]) = 1$ for $1 \leq i \leq n$, then*

$$\text{(5.1)} \qquad P\left(\left|\sum_{i=1}^n X_i\right|\Big/n \geq x\right) \leq 2\exp\{-nx^2(2(b-a)^2)^{-1}\}$$

*for all $x \geq 0$. In particular, when $n = 1$ each $X_i$ is sub-Gaussian with relevant constants $c = 2$ and $k = (2(b-a)^2)^{-1}$. If*

$$\text{(5.2)} \qquad P(|X_i| \geq x) \leq ce^{-kx^2}$$

*for $1 \leq i \leq n$ and all $x \geq 0$, then*

$$\text{(5.3)} \qquad P\left(\left|\sum_{i=1}^n X_i\right|\Big/n \geq x\right) \leq 2\exp\{-nkx^2/(16c)\}.$$

PROOF. First, observe that if $Y$ is a mean zero random variable, then Jensen's inequality implies $E(e^{tY}) \geq e^{tE(Y)} = 1$ for all real $t$. Thus, for $Y_1, Y_2$ independent copies of $Y$, we have $E((Y_1 - Y_2)^l) = 0$ for $l$ odd, and therefore

$$(5.4) \qquad E(e^{tY}) \leq E(e^{t(Y_1 - Y_2)}) = 1 + \sum_{l \geq 1} t^{2l} E((Y_1 - Y_2)^{2l})/(2l)!.$$

If $P(X_i \in [a, b]) = 1$ for $1 \leq i \leq n$, then $E((Y_1 - Y_2)^{2l}) \leq (b - a)^{2l}$ and since $(2l)! \geq 2^l (l!)^2$ for $l \geq 1$ we therefore have

$$E(e^{tY}) \leq 1 + \sum_{l \geq 1} t^{2l} (b - a)^{2l}/(2^l l!) = e^{t^2 (b-a)^2/2}.$$

Applying this estimate to each of the $X_i$'s for $1 \leq i \leq n$, the independence of the $X_i$'s and Markov's inequality implies that for each $t \geq 0$ we have

$$P\left(\sum_{i=1}^n X_i/n \geq x\right) \leq e^{-ntx} \prod_{i=1}^n E(e^{tX_i}) \leq e^{-n(tx - t^2(b-a)^2/2)}.$$

Since $x \geq 0$, minimizing the right-hand side term over $t \geq 0$, we take $t = x/(b-a)^2$, and hence

$$P\left(\sum_{i=1}^n X_i/n \geq x\right) \leq e^{-nx^2/(2(b-a)^2)}.$$

Applying the previous argument to $-\sum_{i=1}^n X_i$, we thus have (5.1).

To prove (5.3), we first show that if $E(Y) = 0$ and

$$(5.5) \qquad\qquad\qquad P(|Y| \geq x) \leq ce^{-kx^2}$$

holds for all $x \geq 0$, then

$$(5.6) \qquad\qquad\qquad E(e^{tY}) \leq e^{4ct^2/k}$$

for all $t \geq 0$. This can be done by utilizing (5.4) and by showing for $Y_1, Y_2$ independent copies of $Y$ that

$$E((Y_1 - Y_2)^{2l}) \leq 2c(4/k)^l l \int_0^\infty e^{-s} s^{l-1} \, ds = 2c(4/k)^l l! \leq (8c/k)^l l!.$$

Thus, (5.6) holds and by applying the previous inequality, independence, and Markov's inequality as before, we have

$$P\left(\sum_{i=1}^n X_i \geq nx\right) \leq \exp\{-nkx^2/(16c)\}.$$

Applying the previous argument to $-\sum_{i=1}^n X_i$, we thus have (5.3), and the lemma is proven. $\square$

LEMMA 5.2. *Let $\{\mathbf{X}_{n,i} : 1 \leq i \leq n\}$ be defined as in (1.2), and assume (2.1) holds for $r = 2$, and the constants $c_{n,j}, k_{n,j}$ are such that $1 \leq c_{n,j}$ and $0 < k_{n,j} \leq \infty$. If $Q_d(\mathbf{x}) = \sum_{j \geq d+1} x_j \mathbf{e}_j$ for $\mathbf{x} \in R^\infty$, and $\tilde{\mathbf{S}}_{n,n}$ is given as in (1.7), then for all $d \geq 0$ and $\delta > 0$*

$$(5.7) \qquad P(\|Q_d(\tilde{\mathbf{S}}_{n,n})\|_\infty \geq \delta) \leq \sum_{j \geq d+1} 2\exp\{-\delta^2 k_{n,j}/(16 c_{n,j})\}.$$

*In addition, if the $V_{n,j}$ are replaced by $n^{1/2}$ in $\tilde{\mathbf{S}}_{n,n}$, then again (5.7) holds.*

PROOF. We first establish (5.7) for general $V_{n,j}$. When the $V_{n,j}$ are replaced by $n^{1/2}$, the result then follows by an immediate simplification of this argument.

If $\theta_{n,i,j} = I(j \leq N_{n,i}) R_{n,i,j}$ as indicated, then $P(\theta_{n,i,j} = 1) = p_{n,j} p$, where $p_{n,j} = P(j \leq N_{n,i})$ for $n \geq 1, j \geq 1$, and for $k = 0, 1, \ldots, n$ we define the events

$$(5.8) \qquad E_{k,n,j} = \bigcup_{I \in \mathcal{I}_{k,n,j}} F_I,$$

where $\mathcal{I}_{k,n,j}$ denotes all subsets $I = \{i_1, \ldots, i_k\}$ of size $k$ in $\{1, \ldots, n\}$ and

$$F_I = \{\theta_{n,i,j} = 1 \text{ for all } i \in I \text{ and } \theta_{n,i,j} = 0 \text{ for } i \in \{1, \ldots, n\} \cap I^c\}.$$

Note that $F_I$ depends on $n$ and $j$, but we suppress that in our notation.

Since $V_{n,j} = \max\{1, \sum_{i=1}^n \theta_{n,i,j}\}$ and $\sum_{i=1}^n \xi_{n,i,j} \theta_{n,i,j} = 0$ on $E_{0,n,j}$, we therefore have for each $\delta > 0, n \geq 1$, and $d \geq 0$ that

$$P(\|Q_d(\tilde{\mathbf{S}}_{n,n})\|_\infty \geq \delta) \leq \sum_{j \geq d+1} \sum_{k=1}^n P\left(\left\{\left|\sum_{i=1}^n \xi_{n,i,j} \theta_{n,i,j}\right| \geq \delta V_{n,j}^{1/2}\right\} \cap E_{k,n,j}\right).$$

Now

$$P\left(\left\{\left|\sum_{i=1}^n \xi_{n,i,j} \theta_{n,i,j}\right| \geq \delta V_{n,j}^{1/2}\right\} \cap E_{k,n,j}\right)$$
$$= \sum_{I \in \mathcal{I}_{k,n,j}} P\left(\left\{\left|\sum_{l=1}^k \xi_{n,i_l,j}\right| \geq \delta k^{1/2}\right\} \cap F_I\right)$$

and letting

$$A_n = \sum_{j \geq d+1} \sum_{k=1}^n P\left(\left\{\left|\sum_{i=1}^n \xi_{n,i,j} \theta_{n,i,j}\right| \geq \delta V_{n,j}^{1/2}\right\} \cap E_{k,n,j}\right),$$

we have by using the independence of the various sequences of random variables involved that

$$A_n = \sum_{j \geq d+1} \sum_{k=1}^{n} \sum_{I \in \mathcal{I}_{k \cdot n,j}} P\left(\left\{\left|\sum_{l=1}^{k} \xi_{n,i_l,j}\right| \geq \delta k^{1/2}\right\}\right) P(F_I, I = \{i_1, \ldots, i_k\})$$

$$\leq 2 \sum_{j \geq d+1} \sum_{k=1}^{n} \exp\{-\delta^2 k_{n,j}/16 c_{n,j}\} P(E_{k,n,j}).$$

Of course, in the previous inequality we are applying (5.3) of Lemma 5.1 to estimate $P(\{|\sum_{l=1}^{k} \xi_{n,i_l,j}| \geq \delta k^{1/2}\})$. Thus, we have (5.7) for general $V_{n,j}$.

When the $V_{n,j}$ are replaced by $n$, the proof is immediate since the random variables $\{\xi_{n,i,j}\theta_{n,i,j}: n \geq 1, i \geq 1, j \geq 1\}$ also satisfy (2.1), and hence one can apply (5.3) immediately to obtain (5.7). Hence, the lemma is proven. □

In order that the probability estimate in the previous lemma be useful $k_{n,j}/c_{n,j}$ must be unbounded as $j$ tends to infinity. Our next task is to see what happens if we remove this assumption, and only ask that this ratio is uniformly bounded below by a strictly positive constant. This is the content of our next lemma, which is a modification of Lemma 5.2.

LEMMA 5.3. *Let $\{\mathbf{X}_{n,i}: 1 \leq i \leq n\}$ be defined as in (1.2), and assume (2.1) holds for $r = 2$, and the constants $c_{n,j}, k_{n,j}$ are such that $1 \leq c_{n,j} \leq c < \infty$ and $0 < k \leq k_{n,j} \leq \infty$. If $\tilde{\mathbf{S}}_{n,n}$ is given as in (1.7), and $N_n^* = \max_{1 \leq i \leq n} N_{n,i}$, then*

$$(5.9) \qquad P(\|\tilde{\mathbf{S}}_{n,n}\|_\infty \geq x) \leq 2E(N_n^*) \exp\left\{-\frac{kx^2}{16c}\right\}.$$

*In addition, if the $V_{n,j}$ are replaced by $n^{1/2}$ in $\tilde{\mathbf{S}}_{n,n}$, then again (5.9) holds.*

REMARK 5.1. If $N_{n,i} = p_n$ for $\{i \geq 1, n \geq 1\}$, then (5.9) immediately implies

$$(5.10) \qquad P(\|\tilde{\mathbf{S}}_{n,n}\|_\infty \geq x) \leq 2p_n \exp\left\{-\frac{kx^2}{16c}\right\}$$

and if the $V_{n,j}^{1/2}$ are replaced by $n^{1/2}$ in $\tilde{\mathbf{S}}_{n,n}$, then again (5.10) holds.

PROOF OF LEMMA 5.3. Following the proof of Lemma 5.2, we observe that if $\theta_{n,i,j} = I(j \leq N_{n,i}) R_{n,i,j}$, then $P(\theta_{n,i,j} = 1) = p_{n,j} p$, where $p_{n,j} = P(j \leq N_{n,i})$ for $n \geq 1, j \geq 1$, and for $m = 0, 1, \ldots, n$ we define the events

$$E_{m,n,j} = \bigcup_{I \in \mathcal{I}_{m,n,j}} F_I,$$

ASYMPTOTIC INFERENCE FOR HIGH-DIMENSIONAL DATA 21

where $\mathcal{I}_{m,n,j}$ denotes all subsets $I = \{i_1, \ldots, i_m\}$ of size $m$ in $\{1, \ldots, n\}$ and

$$F_I = \{\theta_{n,i,j} = 1 \text{ for all } i \in I \text{ and } \theta_{n,i,j} = 0 \text{ for } i \in \{1, \ldots, n\} \cap I^c\}.$$

Recall $V_{n,j} = \max\{1, \sum_{i=1}^n \theta_{n,i,j}\}$ and observe that $\sum_{i=1}^n \xi_{n,i,j}\theta_{n,i,j} = 0$ on $E_{0,n,j}$. Hence, for each $x > 0, n \geq 1$,

$$P(\|\tilde{\mathbf{S}}_{n,n}\|_\infty \geq x)$$
$$\leq \sum_{u \geq 1} \sum_{j=1}^u \sum_{m=1}^n P\left(\left\{\left|\sum_{i=1}^n \xi_{n,i,j}\theta_{n,i,j}\right| \geq xV_{n,j}^{1/2}\right\} \cap E_{m,n,j} \cap \{N_n^* = u\}\right).$$

Now setting $B_{n,m} = \{|\sum_{l=1}^m \xi_{n,i_l,j}| \geq xm^{1/2}\}$

$$P\left(\left\{\left|\sum_{i=1}^n \xi_{n,i,j}\theta_{n,i,j}\right| \geq xV_{n,j}^{1/2}\right\} \cap E_{m,n,j} \cap \{N_n^* = u\}\right)$$
$$= \sum_{I \in \mathcal{I}_{m,n,j}} P(B_{n,m} \cap F_I \cap \{N_n^* = u\})$$

and letting

$$A_n = \sum_{u \geq 1} \sum_{j=1}^u \sum_{m=1}^n P\left(\left\{\left|\sum_{i=1}^n \xi_{n,i,j}\theta_{n,i,j}\right| \geq xV_{n,j}^{1/2}\right\} \cap E_{m,n,j} \cap \{N_n^* = u\}\right),$$

we have by using the independence of the various sequences of random variables involved and (5.3) of Lemma 5.1 that

$$A_n \leq 2 \sum_{u \geq 1} \sum_{j=1}^u \sum_{m=1}^n \exp\{-x^2 k_{n,j}/16 c_{n,j}\} P(E_{m,n,j} \cap \{N_n^* = u\}).$$

Thus, by first summing on $m$ and using $c_{n,j} \leq c$ and $k_{n,j} \geq k$ for all $n, j \geq 1$, we have that

(5.11) $\quad P(\|\tilde{\mathbf{S}}_{n,n}\|_\infty \geq x) \leq 2 \sum_{u \geq 1} \sum_{j=1}^u \exp\{-x^2 k/(16c)\} P(N_n^* = u).$

Hence, this implies (5.9), and when the $V_{n,j}$ are replaced by $n$, the proof follows from the ideas used in the general case. Thus, the lemma is proven. □

Next, we turn to a method which will allow us to handle a broader collection of random variables. Here the $\{\xi_{n,i,j}\}$ satisfy (2.1) with $r \in (0, 2)$, or the less restrictive conditions of polynomial decay given in (2.4). Of course, the results depend on the rate of decay of the tails of the $\{\xi_{n,i,j}\}$, but under a variety of assumptions we are able to obtain further consistency results in this setting. The relevant probability inequalities are obtained in our next lemma, and can be viewed as a substitute for those in Lemma 5.1.

LEMMA 5.4. *For each integer $n \geq 1$ let $X_1, \ldots, X_n$ be independent, mean zero random variables, such that for some $r \in (0, 2)$ we have*

$$P(|X_i| \geq x) \leq c e^{-kx^r} \tag{5.12}$$

*for $1 \leq i \leq n$ and all $x \geq 0$. Then for all $x \geq \sqrt{8} M_{c,k,r}/\sqrt{n}$ and all $s \geq 0$*

$$P\left(\left|\sum_{i=1}^n X_i\right| \geq nx\right) \leq 4\exp\left\{-\frac{nx^2}{32s^2}\right\} + 4cn\exp\left\{-\frac{ks^r}{2^r}\right\}, \tag{5.13}$$

*where $M_{c,k,r}^2 = \int_0^\infty c e^{-kx^{r/2}}\,dx < \infty$. In addition, for all $x \geq \sqrt{8} M_{c,k,r}/\sqrt{n}$ and all $s \geq 1$, we also have*

$$P\left(\left|\sum_{i=1}^n X_i\right| \geq nx\right) \leq 4\exp\left\{-\frac{nkx^2}{128cs^2}\right\} + 4cn\exp\left\{-\frac{ks^r}{2^r}\right\}. \tag{5.14}$$

*Moreover, if for some $c > 0$ and $k > 2$, (5.12) is replaced by*

$$P(|X_i| \geq x) \leq \frac{c}{(1+x)^k} \tag{5.15}$$

*for $1 \leq i \leq n$ and all $x \geq 0$, then for all $x \geq \sqrt{8} M_{c,k}/\sqrt{n}$ and all $s \geq 0$*

$$P\left(\left|\sum_{i=1}^n X_i\right| \geq nx\right) \leq 4\exp\left\{-\frac{nx^2}{32s^2}\right\} + \frac{2^{2+k}cn}{(2+s)^k}, \tag{5.16}$$

*where $M_{c,k}^2 = \int_0^\infty \frac{c}{(1+t^{1/2})^k}\,dt < \infty$.*

REMARK 5.2. If we take $s = n^{1/(2+r)} x^{2/(2+r)}$ then for $x \geq \sqrt{8} M_{c,k,r}/\sqrt{n}$, we have that (5.13) implies

$$P\left(\left|\sum_{i=1}^n X_i\right| \geq nx\right) \leq 4\exp\left\{-\frac{n^{r/(2+r)} x^{2r/(2+r)}}{32}\right\}$$
$$+ 4cn\exp\left\{-\frac{kn^{r/(2+r)} x^{2r/(2+r)}}{2^r}\right\},$$

which makes the exponents on the right of comparable size. Since the proof of (5.13) and (5.14) also implies (5.13) and (5.14) when $r = 2$, it is interesting to note that the previous inequality is not as sharp as that in (5.3) in Lemma 5.1 when $r = 2$.

REMARK 5.3. If the median of each $X_i$ is zero, then (5.13) and (5.14) hold for all $x \geq 0$ and $s \geq 0$. That is, when the medians are zero the key inequality (5.18) below follows directly from (5.8) in [4], page 147, without restrictions on $x$. A similar remark holds for (5.14) provided $s \geq 1$.

PROOF OF LEMMA 5.4. First, we observe that for $r$ fixed, uniformly in $i$, $1 \le i \le n$, (5.12) implies $E(X_i^2) \le M_{c,k,r}^2 < \infty$. Hence if $x \ge \sqrt{8} M_{c,k,r}/\sqrt{n}$ we have by Cheyshev's inequality that $P(|\sum_{i=1}^n X_i| \ge nx/2) \le 1/2$. Now let $Y_1, \ldots, Y_n$ be an independent copy of $X_1, \ldots, X_n$ and observe that

$$(5.17) \quad P\left(\left|\sum_{i=1}^n X_i\right| \ge nx\right) P\left(\left|\sum_{i=1}^n Y_i\right| \le nx/2\right) \le P\left(\left|\sum_{i=1}^n (X_i - Y_i)\right| \ge nx/2\right).$$

Then for all $x \ge \sqrt{8} M_{c,k,r}/\sqrt{n}$, we have

$$(5.18) \quad P\left(\left|\sum_{i=1}^n X_i\right| \ge nx\right) \le 2P\left(\left|\sum_{i=1}^n (X_i - Y_i)\right| \ge nx/2\right).$$

Taking $s \ge 0$, we define $(X_i - Y_i)^s = (X_i - Y_i) I(|X_i - Y_i| \le s)$. Then

$$(5.19) \quad P\left(\left|\sum_{i=1}^n (X_i - Y_i)\right| \ge nx/2\right) \le I_n(s,x) + II_n(s,x),$$

where

$$I_n(s,x) = P\left(\left|\sum_{i=1}^n (X_i - Y_i)^s\right| \ge nx/2\right)$$

and

$$II_n(s,x) = P\left(\max_{1 \le i \le n} |(X_i - Y_i) - (X_i - Y_i)^s| > 0\right).$$

Applying (5.1) to $(X_1 - Y_1)^s, \ldots, (X_n - Y_n)^s$, we see that

$$(5.20) \quad I_n(s,x) \le 2\exp\{-n(x/2)^2 (2(2s)^2)^{-1}\}$$

and (5.12) implies

$$(5.21) \quad II_n(s,x) \le 2 \sum_{i=1}^n P(|X_i| > s/2) \le 2cn \exp\left\{-\frac{ks^r}{2^r}\right\}.$$

Applying (5.18), (5.19), (5.20) and (5.21) we thus have (5.13).

The proof of (5.14) follows that for (5.13) up to the point we apply (5.1) of Lemma 5.1 to $I_n(s,x)$ in (5.20). At this point, we now apply (5.3) of Lemma 5.1 to the random variables $(X_1 - Y_1)^s, \ldots, (X_n - Y_n)^s$. That is, (5.12) implies that for all $x \ge 0$ and $1 \le i \le n$ that

$$P(|(X_i - Y_i)^s| \ge x) \le P(|X_i^s| \ge x/2) + P(|Y_i^s| \ge x/2) \le 2ce^{-kx^2/(4s^2)},$$

where the last inequality follows since $x^r/2^r \ge x^2/(4s^2)$ when $0 \le x \le 2s$, $0 < r < 2$, and $s \ge 1$. Hence, by (5.3) of Lemma 5.1, with $k$ replaced by $k/(4s^2)$ and $c$ by $2c$, we obtain

$$(5.22) \quad I_n(s,x) \le 2\exp\{-nkx^2/(128cs^2)\}.$$

Now combining (5.22) and the estimate for $II_n(s,x)$ in (5.21) to (5.18) and (5.19), we obtain (5.14).

Next, we observe that uniformly in i, $1 \le i \le n$, (5.15) and $k > 2$ implies

$$E(X_i^2) = \int_0^\infty P(|X_i| > t^{1/2})\,dt \le \int_0^\infty \frac{c}{(1+t^{1/2})^k}\,dt < \infty.$$

Hence, if we choose $x \ge \sqrt{8}M_{c,k}/\sqrt{n}$ we have by the argument leading to (5.17)–(5.19) that

$$(5.23) \quad P\left(\left|\sum_{i=1}^n (X_i - Y_i)\right| \ge nx/2\right) \le I_n(s,x) + II_n(s,x),$$

where

$$I_n(s,x) = P\left(\left|\sum_{i=1}^n (X_i - Y_i)^s\right| \ge nx/2\right)$$

and

$$II_n(s,x) = P\left(\max_{1\le i\le n} |(X_i - Y_i) - (X_i - Y_i)^s| > 0\right).$$

Applying (5.1) to $(X_1 - Y_1)^s, \ldots, (X_n - Y_n)^s$, we see that

$$(5.24) \quad I_n(s,x) \le 2\exp\{-n(x/2)^2(2(2s)^2)^{-1}\}$$

and (5.14) implies

$$(5.25) \quad II_n(s,x) \le \sum_{i=1}^n P(|X_i - Y_i| > s) \le 2\sum_{i=1}^n P(|X_i| > s/2) \le \frac{2cn}{(1+s/2)^k}.$$

Applying (5.23), (5.24) and (5.25), we have (5.16). Thus, Lemma 5.4 is proven. □

## 6. Proofs of consistency results.

6.1. *Proof of Theorem 2.1.* Applying Lemma 5.2 with $d = 0$, we have for all $x > 0$ and each integer $n \ge 1$ that

$$(6.1) \quad P(\|\tilde{\mathbf{S}}_{n,n}\|_\infty \ge x) \le \sum_{j\ge 1} 2\exp\{-x^2 k_{n,j}/(16c_{n,j})\}.$$

Taking $x = \varepsilon a_n$ in (6.1) and applying (2.5), we thus have (2.6) for general $V_{n,j}^{1/2}$, and also when the $V_{n,j}^{1/2}$ are replaced by $n^{1/2}$.

If the constants $c_{n,j}$ and $k_{n,j}$ satisfy (2.7) as indicated, then with $a_n = (L(n+3))^{1/2}$ and $x = \varepsilon a_n$ in (6.1) we have

$$P(\|\tilde{\mathbf{S}}_{n,n}\|_\infty \geq \varepsilon(L(n+3))^{1/2}) \leq \sum_{j\geq 1} 2\exp\{-\varepsilon^2 \delta L(n+3)L(j+3)\}$$
(6.2)
$$= 2\sum_{j\geq 1}(j+3)^{-\varepsilon^2 \delta L(n+3)}.$$

Thus, for $\varepsilon^2 \delta > 1$,

$$\sum_{n\geq 1} P(\|\tilde{\mathbf{S}}_{n,n}\|_\infty \geq \varepsilon(L(n+3))^{1/2})$$
(6.3)
$$\leq 2\sum_{n\geq 1}\int_3^\infty x^{-\varepsilon^2 \delta L(n+3)}\,dx$$

$$\leq 2\sum_{n\geq 1}\frac{3^{-L(n+3)-1}}{(L(n+3)-1)} < \infty$$

and hence (2.8) holds for general $V_{n,j}$. In particular, we then have from (6.3) that (2.9) is immediate, and it remains to show $E(e^{\alpha M^2}) < \infty$ for all $\alpha > 0$ sufficiently small. Now

$$E(e^{\alpha M^2}) = \int_0^\infty P(e^{\alpha M^2} > t)\,dt \leq 3 + \int_3^\infty P\left(M > \left(\frac{\log t}{\alpha}\right)^{1/2}\right)dt$$

and

$$\int_3^\infty P\left(M > \left(\frac{\log t}{\alpha}\right)^{1/2}\right)dt$$

$$\leq \sum_{n\geq 1}\int_3^\infty P\left(\|\tilde{\mathbf{S}}_{n,n}\|_\infty \geq (L(n+3))^{1/2}\left(\frac{\log t}{\alpha}\right)^{1/2}\right)dt$$

$$\leq 2\sum_{n\geq 1}\sum_{j\geq 1}\int_3^\infty \exp\left\{-\delta L(j+3)L(n+3)\frac{\log t}{\alpha}\right\}dt,$$

where the last inequality follows from (6.1) and that (2.7) holds.

Therefore, for $\alpha < \delta/2$ we have

$$E(e^{\alpha M^2}) \leq 3 + 2\sum_{n\geq 1}\sum_{j\geq 1}\int_3^\infty \exp\{-2L(j+3)L(n+3)\log t\}\,dt$$

$$= 3 + 6\sum_{n\geq 1}\sum_{j\geq 1}\frac{\exp\{-2\log 3 L(j+3)L(n+3)\}}{2L(j+3)L(n+3)-1}.$$

Now $x, y \geq 1 + \eta$ for some $\eta > 0$ implies $xy \geq (x+y)(1+\eta)/2$ and hence since $j, n \geq 1$ implies $L(j+3), L(n+3) \geq L4 \geq 1 + \eta$ for $\eta = L4 - 1 > 0$ we have

$$E(e^{\alpha M^2}) \leq 3 + 6 \sum_{n \geq 1} \sum_{j \geq 1} \frac{\exp\{-(1+\eta)\log 3(L(j+3) + L(n+3))\}}{2L(j+3)L(n+3) - 1} < \infty$$

since $(1+\eta)\log 3 > 1$. Since (5.1) holds when the $V_{n,j}^{1/2}$ are replaced by $n^{1/2}$, the proof also holds in this situation. Thus, Theorem 2.1 is proven.

6.2. *Proof of Theorem 2.2.* Under our assumptions, Lemma 5.3 with $x = h(n)$ implies that

$$P(\|\tilde{\mathbf{S}}_{n,n}\|_\infty \geq h(n)) \leq 2E(N_n^*) \exp\left\{-\frac{kh(n)^2}{16c}\right\}.$$

Since $x = h(n) = (\theta_1^{-1} L(E(N_n^*)) + \theta_2 Ln)^{1/2}$, $\theta_1 = k/(16c)$ and $\theta_2 > 0$, we have

$$P(\|\tilde{\mathbf{S}}_{n,n}\|_\infty \geq h(n)) \leq 2L(E(N_n^*)) \exp\{-L(E(N_n^*)) - \theta_1 \theta_2 Ln\}.$$

Since $\theta_1 > 0$, we thus have (2.13) if $\theta_2 > 0$, and (2.14) follows immediately provided $\theta_1 \theta_2 > 1$. Since the above holds for general $V_{n,j}^{1/2}$, and also the $n^{1/2}$ normalizations, Theorem 2.2 is proven.

6.3. *Proof of Theorem 2.3.* First, observe that for all $x \geq 0$ that

$$(6.4) \qquad P(\|\mathbf{S}_{n,n}\|_\infty \geq nx) = P\left(\bigcup_{j=1}^{N_n^*}\left\{\left|\sum_{i=1}^n \xi_{n,i,j}\theta_{n,i,j}\right| \geq nx\right\}\right).$$

Letting $\mathbf{b} = (b_1, \ldots, b_n)$, where $b_i$ is a positive integer for $1 \leq i \leq n$, and setting $E_{n,\mathbf{b}} = \{N_{n,1} = b_1, \ldots, N_{n,n} = b_n\}$, we thus have by conditioning on $E_{n,\mathbf{b}}$ that

$$P(\|\mathbf{S}_{n,n}\|_\infty \geq nx) \leq \sum_{(b_1,\ldots,b_n)} \sum_{j=1}^{\max(b_1,\ldots,b_n)} J(n, j, \mathbf{b}, x) P(E_{n,\mathbf{b}}),$$

where $J(n, j, \mathbf{b}, x) = P(|\sum_{i=1}^n \xi_{n,i,j}\theta_{n,i,j}| \geq nx | E_{n,\mathbf{b}})$. Fixing $n$ and $j$, and defining $X_i = \xi_{n,i,j}\theta_{n,i,j} I(j \leq b_i)$ for $1 \leq i \leq n$, we see $X_1, \ldots, X_n$ are independent random variables, and it is easy to check from our assumptions on $c_{n,j}$ and $k_{n,j}$, and (2.1), that for all $x \geq 0$

$$P(|X_i| \geq x) \leq ce^{-kx^r}.$$

Therefore, $X_1, \ldots, X_n$ satisfy the conditions in Lemma 5.4 and using the independence of the sequences $\{\xi_{n,i,j}\}$, $\{R_{n,i,j}\}$, and $\{N_{n,i}\}$, we have for $x \geq \sqrt{8} M_{c,k,r}/\sqrt{n}$ and $s \geq 1$ that (5.14) implies

$$J(n, j, \mathbf{b}, x) \leq 4 \exp\{-nkx^2/(128cs^2)\} + 4cn \exp\left\{-\frac{ks^r}{2^r}\right\}.$$

Combining the previous inequalities in this proof, we have that

(6.5) $\qquad P(\|\mathbf{S}_{n,n}\|_\infty \geq nx) \leq 4E(N_n^*)[A_n(1) + A_n(2)],$

where $A_n(1) = \exp\{-nkx^2/(128cs^2)\}$ and $A_n(2) = cn \exp\{-\frac{ks^r}{2^r}\}$. Recalling $h(n) = (c_2^{-1} L(E(N_n^*)) + c_3 L(n))^{1/2}$ and taking $s = s_n = c_1(L(E(N_n^*)) + 2 \times L(n))^{1/r}$ and $x = x_n = s_n \{c_2^{-1} L(E(N_n^*)) + c_3 L(n)\}^{1/2}/n^{1/2}$ in (6.5), then for all sufficiently large $n$ we have $x \geq \sqrt{8} M_{c,k,r}/\sqrt{n}$, $s \geq 1$, and (6.5) holds. Thus, (2.18) holds if $c_1 > 2/k^{1/r}$ and $c_3 > 0$, and (2.19) follows if we also have $c_2 c_3 > 1$. Thus, Theorem 2.3 is proven when (2.1) holds and $0 < r < 2$.

We now turn to the situation where there is only polynomial decay in the tails of the data $\xi_{n,i,j}$ as in (2.4), where $c_{n,j} \leq c$ and $k_{n,j} \geq k$ for all $n \geq 1, j \geq 1$ and $1 \leq c < \infty$, $2 < k < \infty$. Then the random variables $\xi_{n,i,j} \theta_{n,i,j}$ are also easily seen to satisfy (2.4), and arguing as in (6.4) and (6.5), and applying (5.16), we have for $s \geq 0$ and $x \geq \sqrt{8} M_{c,k}/n^{1/2}$ that

(6.6) $\qquad P(\|\mathbf{S}_{n,n}\|_\infty \geq nx) \leq 4E(N_n^*)\left[\exp\left\{-\frac{nx^2}{32s^2}\right\} + \frac{2^k cn}{(2+s)^k}\right].$

Taking $s = s_n = (nE(N_n^*))^{1/k+\beta}$, $\beta > 0$, and $x = x_n = bs_n(LE(N_n^*))^{1/2}/n^{1/2}$, then for all $n$ sufficiently large

$$P(\|\mathbf{S}_{n,n}\|_\infty \geq bn^{1/2} s_n (LE(N_n^*))^{1/2}) \leq 4E(N_n^*)[A_n(3) + A_n(4)],$$

where $A_n(3) = \exp\{-\frac{b^2}{32} L(E(N_n^*))\}$ and $A_n(4) = 2^k cn[2+s_n]^{-k}$. Thus, (2.20) holds when $\beta > 0$ by taking $b$ large and using the definition of $s_n$. Moreover, (2.21) holds if $b \geq 8$ and $k\beta > 1/2$, and (2.22) holds when if $b > 8$ and $\beta k(\gamma + 1) > 1$. Thus, Theorem 2.3 is proven.

**7. Proof of Theorem 2.4.** The proof proceeds with a sequence of lemmas.

Our first lemma provides tightness, and shows $\tilde{\mathbf{T}}_{n,n}$ converges in probability to zero in $c_0$.

LEMMA 7.1. *Under the conditions (2.29) and (2.30),*

(7.1) $\qquad \{\mathcal{L}(\tilde{\mathbf{S}}_{n,n}) : n \geq 1\}$ *is tight in* $c_0$.

*In addition, if for each $j \geq 1$ we have $\lim_{n \to \infty} P(N_{n,1} < j) = 0$, then $\tilde{\mathbf{T}}_{n,n}$ converges in probability to zero in $c_0$.*

PROOF. For general $V_{n,j}$, or if the $V_{n,j}$ are replaced by $n$, (5.7) implies that

$$(7.2) \qquad P(\|Q_d(\tilde{\mathbf{S}}_{n,n})\|_\infty \geq \delta) \leq \sum_{j \geq d+1} 2\exp\{-\delta^2 k_{n,j}/(16 c_{n,j})\}.$$

Hence, (2.30) implies for $\delta > 0$ arbitrary that

$$(7.3) \qquad \lim_{d \to \infty} \sup_{n \geq 1} P(\|Q_d(\tilde{\mathbf{S}}_{n,n})\|_\infty > \delta) = 0.$$

Now (2.1) easily implies $E(\xi_{n,i,j}^2) \leq c_{n,j}/k_{n,j}$, and the independence of the sequences of random variables involved implies for each $j \geq 1$ that

$$P\left(\left|\sum_{i=1}^n \xi_{n,i,j}\theta_{n,i,j}\right| \geq bV_{n,j}^{1/2}\right) \leq b^{-2} c_{n,j}/k_{n,j}.$$

Thus (2.29), (7.3), and an application of the remark on page 49 of [13] easily combine to prove the tightness in (7.1) for general $V_{n,j}$ and also when the $V_{n,j}$ are replaced by $n$.

If $\tilde{\mathbf{T}}_{n,n}$ is defined as in (1.9), then for each $\varepsilon > 0$ and $d \geq 1$ we have

$$P(\|\tilde{\mathbf{T}}_{n,n}\|_\infty > 2\varepsilon) \leq P(\|\tilde{\mathbf{T}}_{n,n} - Q_d(\tilde{\mathbf{T}}_{n,n})\|_\infty > \varepsilon) + P(\|Q_d(\tilde{\mathbf{T}}_{n,n})\|_\infty > \varepsilon).$$

Now (7.3) immediately implies there exists $d \geq 1$ such that uniformly in $n$, $P(\|Q_d(\tilde{\mathbf{T}}_{n,n})\|_\infty > \varepsilon) < \varepsilon/2$. Independence of the sequences of random variables involved then implies for each $j \geq 1$ and $b > 0$ that

$$P\left(\left|\sum_{i=1}^n \xi_{n,i,j}\theta_{n,i,j}\right| \geq bV_{n,j}\right)$$
$$\leq b^{-2} E\left(E\left(\left(\sum_{i=1}^n \xi_{n,i,j}\theta_{n,i,j}\right)^2 \bigg/ V_{n,j}^2 \bigg| \theta_{n,1,j},\ldots,\theta_{n,n,j}\right)\right)$$
$$\leq b^{-2} E(V_{n,j}^{-1}) c_{n,j}/k_{n,j}.$$

Since $\lim_{n \to \infty} P(N_{n,1} < j) = 0$, the weak law of large numbers and Chebyshev's inequality applied to the i.i.d. sequence of random variables $\{\theta_{n,i,j} : i \geq 1\}$ implies for each fixed $j \geq 1$ and $M > 0$ that $\limsup_{n \to \infty} P(V_{n,j} \leq M) = 0$. Thus, $\limsup_{n \to \infty} E(V_{n,j}^{-1}) = 0$, so for each fixed $j \geq 1$ and $b > 0$ we have

$$\lim_{n \to \infty} P\left(\left|\sum_{i=1}^n \xi_{n,i,j}\theta_{n,i,j}\right| \geq bV_{n,j}\right) = 0.$$

Now

$$P(\|\tilde{\mathbf{T}}_{n,n} - Q_d(\tilde{\mathbf{T}}_{n,n})\|_\infty > \varepsilon) \leq \sum_{j=1}^d P\left(\left|\sum_{i=1}^n \xi_{n,i,j}\theta_{n,i,j}\right| \geq \varepsilon V_{n,j}\right)$$

and hence from the above we have for each $\varepsilon > 0$ that

$$\lim_{n \to \infty} P(\|\tilde{\mathbf{T}}_{n,n}\|_\infty > 2\varepsilon) \leq \varepsilon.$$

Thus, the lemma is proven. □

Now that we have tightness of $\{\mathcal{L}(\tilde{\mathbf{S}}_{n,n}) : n \geq 1\}$ in $c_0$, the next step of the proof is to show that the finite-dimensional distributions induced by $\bigcup_{d \geq 1} c_{0,d}^*$ are the same for every limiting measure of $\{\mathcal{L}(\tilde{\mathbf{S}}_{n,n}) : n \geq 1\}$. Here $c_0^*$ denotes the continuous linear functionals on $c_0$ and

(7.4) $\qquad c_{0,d}^* = \{f \in c_0^* : f(Q_d(\mathbf{x})) = 0 \text{ for all } \mathbf{x} \in c_0\}.$

We start by showing that the limiting covariance functions $\Gamma(k, \cdot, \cdot)$ given in (2.26)–(2.28) determine the limiting variance of $f(\tilde{\mathbf{S}}_{n,n})$ for each $d \geq 1$ and $f \in c_{0,d}^*$. This follows from our next lemma.

LEMMA 7.2. *If (2.26)–(2.28) hold and $P(\min_{1 \leq i \leq n} N_{n,i} < d) = o(1/n^2)$ as $n$ tends to infinity, then for all $d \geq 1$ and $f \in c_{0,d}^*$ we have*

(7.5) $\qquad \lim_{n \to \infty} E(f^2(\tilde{\mathbf{S}}_{n,n})) = \sum_{u=1}^{d} \sum_{v=1}^{d} \Gamma(1, u, v) f(\mathbf{e}_u) f(\mathbf{e}_v).$

*If $V_{n,j}$ is replaced by $n$ in $\tilde{\mathbf{S}}_{n,n}$, then (7.5) holds with $\Gamma(1, \cdot, \cdot)$ replaced by $\Gamma(2, \cdot, \cdot)$ in the right-hand side term.*

PROOF. Since $f \in c_{0,d}^*$, we have

(7.6) $\quad E(f^2(\tilde{\mathbf{S}}_{n,n})) = \sum_{u=1}^{d} \sum_{v=1}^{d} \left[ \sum_{i=1}^{n} \Gamma_{n,i}(u, v) E\left( \frac{\theta_{n,i,u}}{V_{n,u}^{1/2}} \frac{\theta_{n,i,v}}{V_{n,v}^{1/2}} \right) \right] f(\mathbf{e}_u) f(\mathbf{e}_v),$

where $\Gamma_{n,i}(u, v) = E(\xi_{n,i,u} \xi_{n,i,v})$. This follows immediately since the sequences $\{\theta_{n,i,j}\}$ and $\{V_{n,j}\}$ are independent of the sequence $\{\xi_{n,i,j}\}$, and $E(\xi_{n,i,j}) = 0$ with the random variables $\xi_{n,i,j}$ independent in $i$. Hence, (7.5) and Lemma 7.1 follows from (7.6) once we prove the following lemma. The situation when $V_{n,j}$ is replaced by $n$ in $\tilde{\mathbf{S}}_{n,n}$ is simpler, so for the time being, we assume the $V_{n,j}$ are random. The nonrandom case will be taken up later. □

LEMMA 7.3. *Under the assumptions of Theorem 2.4, we have*

(7.7) $\qquad \lim_{n \to \infty} \sum_{i=1}^{n} \Gamma_{n,i}(u, v) E\left( \frac{\theta_{n,i,u}}{V_{n,u}^{1/2}} \frac{\theta_{n,i,v}}{V_{n,v}^{1/2}} \right) = \Gamma(1, u, v).$

PROOF. Since $\lim_{n\to\infty} \sum_{i=1}^{n} \Gamma_{n,i}(u,v)/n = \Gamma(u,v)$ by assumption, (7.7) will follow if we first show for $u \neq v$ that

(7.8) $$E\left(\frac{\theta_{n,i,u}}{V_{n,u}^{1/2}} \frac{\theta_{n,i,v}}{V_{n,v}^{1/2}}\right) = \frac{p}{n} a_{n,i},$$

where $\lim_{n\to\infty} \sup_{1\leq i \leq n} |a_{n,i} - 1| = 0$, and that

(7.9) $$\sup_{n\geq 1, i\geq 1} \Gamma_{n,i}(u,v) < \infty.$$

Now

(7.10) $$\sup_{n\geq 1, i\geq 1} \Gamma_{n,i}(u,v) = \sup_{n\geq 1, i\geq 1} E(\xi_{n,i,u}\xi_{n,i,v})$$

and, as mentioned earlier, (2.1) with $r = 2$ implies $E(\xi_{n,i,u}^2) \leq c_{n,j}/k_{n,j}$. Therefore, the Cauchy–Schwarz inequality and (2.29) easily combine to imply (7.9). Hence, when $u \neq v$ it remains to prove (7.8).

To verify (7.8) for random $V_{n,j}$, for $n \geq 1, u \geq 1$, let

$$\Lambda_{n,i} = \frac{\theta_{n,i,u}}{V_{n,u}^{1/2}} \frac{\theta_{n,i,v}}{V_{n,v}^{1/2}} \quad \text{and} \quad W_{n,u} = \max\left\{1, \sum_{i=1}^{n} R_{n,i,u}\right\}.$$

Thus, we have

(7.11) $$E(\Lambda_{n,i}) = E\left(\Lambda_{n,i} I\left(\min_{1\leq i\leq n} N_{n,i} \geq d\right)\right) + E\left(\Lambda_{n,i} I\left(\min_{1\leq i\leq n} N_{n,i} < d\right)\right)$$
$$= B_n(1,i) - B_n(2,i) + B_n(3,i),$$

where

$$B_n(1,i) = E\left(\frac{R_{n,i,u}}{W_{n,u}^{1/2}} \frac{R_{n,i,v}}{W_{n,v}^{1/2}}\right), \qquad B_n(3,i) = E\left(\Lambda_{n,i} I\left(\min_{1\leq i\leq n} N_{n,i} < d\right)\right)$$

and

$$B_n(2,i) = E\left(\frac{R_{n,i,u}}{W_{n,u}^{1/2}} \frac{R_{n,i,v}}{W_{n,v}^{1/2}} I\left(\min_{1\leq i\leq n} N_{n,i} < d\right)\right).$$

Also,

(7.12) $$B_n(3,i) \leq P\left(\min_{1\leq i\leq n} N_{n,i} < d\right) = o(1/n)$$

and

(7.13) $$B_n(2,i) \leq P\left(\min_{1\leq i\leq n} N_{n,i} < d\right) = o(1/n).$$

Hence, (7.8) will follow if we show for all $i, 1 \leq i \leq n$, and $u \neq v$ that

(7.14) $$B_n(1,i) = \frac{p}{n} a_n,$$

where $\lim_{\to\infty} a_n = 1$. To verify (7.14), we establish the following lemma, which immediately implies (7.14) when $u \neq v$. If $u = v$, then from the above we see that the analogue of (7.14) required is that $E(R_{n,i,u}/W_{n,u}) = c_n/n$ where $\lim_{n\to\infty} c_n = 1$. This follows since from the proof of Lemma 7.4 below we actually have $c_n = 1 - (1-p)^n$. Hence, the proof of Lemma 7.3 and also Lemma 7.2 for random $V_{n,j}$, will follow once Lemma 7.4 is established. □

LEMMA 7.4. *Let $\{X_i : 1 \leq i \leq n\}$ and $\{Y_i : 1 \leq i \leq n\}$ be independent collections of i.i.d. Bernoulii random variables with $P(X_i = 1) = P(Y_i = 1) = p$. Let $A_n = \max\{1, \sum_{i=1}^n X_i\}$ and $B_n = \max\{1, \sum_{i=1}^n Y_i\}$. Then for all $i, 1 \leq i \leq n$, we have*

$$E\left(\frac{X_i}{A_n^{1/2}} \frac{Y_i}{B_n^{1/2}}\right) = \frac{p}{n} b_n, \tag{7.15}$$

*where $\lim_{\to\infty} b_n = 1$.*

PROOF. By the independence assumed, we have

$$E\left(\frac{X_i}{A_n^{1/2}} \frac{Y_i}{B_n^{1/2}}\right) = E\left(\frac{X_i}{A_n^{1/2}}\right) E\left(\frac{Y_i}{B_n^{1/2}}\right)$$

and hence since $\{X_i : 1 \leq i \leq n\}$ and $\{Y_i : 1 \leq i \leq n\}$ i.i.d. Bernoulii random variables with $P(X_i = 1) = P(Y_i = 1) = p$, it suffices to verify that

$$E\left(\frac{X_1}{A_n^{1/2}}\right) = \left(\frac{p}{n}\right)^{1/2} c_n, \tag{7.16}$$

where $\lim_{\to\infty} c_n = 1$. Now, using Jensen's inequality and an easy calculation, we have

$$E\left(\frac{X_1}{A_n^{1/2}}\right) \leq p\left(E\left(\frac{1}{1+\sum_{i=2}^n X_i}\right)\right)^{1/2} = p\left(\frac{1-(1-p)^n}{np}\right)^{1/2}.$$

Hence, since we assume $0 < p \leq 1$, we have

$$E\left(\frac{X_i}{A_n^{1/2}}\right) \leq \left(\frac{p}{n}\right)^{1/2}. \tag{7.17}$$

Thus, Lemma 7.4 will follow provided we establish a comparable lower bound.

Now for each $\varepsilon, 0 < \varepsilon < p$, we have

$$E\left(\frac{X_1}{A_n^{1/2}}\right) \geq \frac{p}{C(n,p,\varepsilon)} \sum_{\{k \,:\, |k/(n-1)-p| < \varepsilon/L(n)\}} P\left(\sum_{i=2}^n X_i = k\right),$$

32    J. KUELBS AND A. N. VIDYASHANKARwhere $C(n,p,\varepsilon) = (1 + (n-1)p - (n-1)\varepsilon/L(n))^{1/2}$. Since $\varepsilon > 0$, Theorem 1 of [12] implies
$$\sum_{\{k:\,|k/(n-1)-p|<\varepsilon/L(n)\}} P\left(\sum_{i=2}^{n} X_i = k\right) \geq 1 - 2\exp\left\{-2(n-1)\left(\frac{\varepsilon}{L(n)}\right)^2\right\}$$

and hence

$$E\left(\frac{X_1}{A_n^{1/2}}\right) \geq \frac{p}{C(n,p,\varepsilon)}\left[1 - 2\exp\left\{-2(n-1)\left(\frac{\varepsilon}{L(n)}\right)^2\right\}\right] = \left(\frac{p}{n}\right)^{1/2} d_n,$$

where $\lim_{n\to\infty} d_n = 1$. This implies the comparable lower bound to (7.17), and therefore Lemma 7.4 holds.

As mentioned earlier, Lemma 7.4 completes the proof of Lemma 7.3, and hence Lemma 7.2 is established with (7.5) providing a limiting variance function when the $V_{n,j}$ are random. If the $V_{n,j}$ are replaced by $n$ in $\tilde{\mathbf{S}}_{n,n}$, then the proof of Lemma 7.3 with the right-hand side of (7.7) being $\Gamma(2,u,v)$ is much simpler, and the details are left for the reader. Hence, Lemma 7.2 is proven. $\square$

Now that Lemma 7.2 is verified, the next step is to show for all $d \geq 1$, $f \in c_{0,d}^*$, and random $V_{n,j}$ that all limit laws of $\{\mathcal{L}(f(\tilde{\mathbf{S}}_{n,n})): n \geq 1\}$ are centered Gaussian random variables with variance given by

(7.18) $$\sigma^2(f) = \sum_{u=1}^{d}\sum_{v=1}^{d} \Gamma(1,u,v) f(\mathbf{e}_u) f(\mathbf{e}_v).$$

Of course, if the $V_{n,j}$ are replaced by $n$ in $\tilde{\mathbf{S}}_{n,n}$, then (7.18) holds with $\Gamma(1,u,v)$ replaced by $\Gamma(2,u,v)$.

To verify this step of the proof, we first prove a lemma which will put us in position to allow an application of Lyapunov's central limit theorem.

LEMMA 7.5.  *For each integer $d \geq 1$ and $\mathbf{x} \in c_0$, let $\Pi_d(\mathbf{x}) = \sum_{j=1}^{d} x_j \mathbf{e}_j$. Under the conditions of the theorem, we have for each $d \geq 1$ that*

(7.19) $$\lim_{n\to\infty} \sum_{i=1}^{n} E(\|\Pi_d(\tilde{\mathbf{X}}_{n,i})\|_\infty^4) = 0,$$

*where*

$$\tilde{\mathbf{X}}_{n,i} = \sum_{j\geq 1} \frac{\xi_{n,i,j}\theta_{n,i,j}}{V_{n,j}^{1/2}} \mathbf{e}_j.$$

PROOF. By Jensen's inequality, we see that

$$\|\Pi_d(\tilde{\mathbf{X}}_{n,i})\|_\infty^4 \leq \left|\sum_{j=1}^d \frac{|\xi_{n,i,j}|\theta_{n,i,j}}{V_{n,j}^{1/2}}\right|^4 \leq 2^{3(d-1)} \sum_{j=1}^d \frac{\xi_{n,i,j}^4 \theta_{n,i,j}}{V_{n,j}^2}.$$

Hence,

$$E(\|\Pi_d(\tilde{\mathbf{X}}_{n,i})\|_\infty^4) \leq 2^{3(d-1)} \sum_{j=1}^d E\left(\frac{\xi_{n,i,j}^4 \theta_{n,i,j}}{V_{n,j}^2}\right)$$

and the lemma will follow if we show

(7.20) $$\lim_{n\to\infty} \sum_{i=1}^n E\left(\frac{\xi_{n,i,j}^4 \theta_{n,i,j}}{V_{n,j}^2}\right) = 0$$

for $j = 1, \ldots, d$ and all $d \geq 1$. Now $E(\frac{\xi_{n,i,j}^4 \theta_{n,i,j}}{V_{n,j}^2}) \leq (E(\xi_{n,i,j}^8))^{1/2}(E(\frac{\theta_{n,i,j}}{V_{n,j}^4}))^{1/2}$, and using (2.1) with $r = 2$ we have $E(\xi_{n,i,j}^8) \leq 24 c_{n,j}/k_{n,j}^4$. Applying (2.29) and that $c_{n,j} \geq 1$ for all $n \geq 1, j \geq 1$ we therefore have, uniformly in $n$ and $j$, that $E(\xi_{n,i,j}^8) < \infty$. Moreover, since $P(\min_{1\leq i \leq n} N_{n,i} < d) = o(1/n^2)$, one can show that

$$E\left(\frac{\theta_{n,i,j}}{V_{n,j}^4}\right) = pE\left(\frac{1}{(1 + \sum_{i=2}^n X_i)^4}\right) + o(1/n^2).$$

Therefore, let

$$A_n = \sum_{k=0}^n \frac{1}{(1+k)^4} P\left(\sum_{i=1}^n X_i = k\right)$$

and we want an appropriate upper bound on $A_{n-1}$. Now $A_n \leq B_n + 2\exp\{-2n/(L(n))^2\}$, where

$$B_n = \sum_{\{k : |k/n - p| \leq 1/L(n)\}} \frac{1}{(1+k)^4} P\left(\sum_{i=1}^n X_i = k\right)$$

and the exponential term follows from an immediate application of Theorem 1 in [12]. Now

$$B_n \leq \frac{1}{(np)^4[1 + 1/(np) - 1/(pL(n))]^4} \leq \frac{2}{(np)^4}$$

for all $n$ sufficiently large. Therefore, for all $n \geq n_0$, we have

$$A_{n-1} \leq \frac{2}{((n-1)p)^4} + 2\exp\left\{-\frac{2(n-1)}{(L(n-1))^2}\right\},$$

which implies

$$E\left(\frac{\theta_{n,i,j}}{V_{n,j}^4}\right) = E\left(\frac{I(j \leq N_{n,i})R_{n,i,j}}{V_{n,j}^4}\right) \leq \frac{2}{((n-1)p)^4} + o(1/n^2)$$

uniformly in $i \geq 1, j \geq 1$. Thus, uniformly in $i,j \geq 1$, we have

$$E\left(\frac{\xi_{n,i,j}^4 \theta_{n,i,j}}{V_{n,j}^2}\right) \leq (E(\xi_{n,i,j}^8))^{1/2}\left(E\left(\frac{\theta_{n,i,j}}{V_{n,j}^4}\right)\right)^{1/2} = o(1/n),$$

which implies (7.20). Thus, (7.19) holds by the inequality prior to (7.20), and Lemma 7.5 is proven for random $V_{n,j}$. If the $V_{n,j}$ are replaced by $n$ in $\tilde{\mathbf{S}}_{n,n}$, then the proof is even easier and details are left to the reader. Hence, Lemma 7.5 holds for both normalizations. □

The next lemma completes the proof of Theorem 2.4.

LEMMA 7.6. *The functions $\Gamma(1,\cdot,\cdot)$ and $\Gamma(2,\cdot,\cdot)$ defined by (2.26)–(2.28), are covariances of centered Gaussian measures $\gamma_1$ and $\gamma_2$, respectively, on $c_0$. Furthermore, if the $V_{n,j}$ are random, then $\tilde{\mathbf{S}}_{n,n}$ converges weakly to $\gamma_1$ on $c_0$, and if the $V_{n,j}$ are replaced by $n$, then $\tilde{\mathbf{S}}_{n,n}$ converges weakly to $\gamma_2$ on $c_0$. In addition, for each $f \in c_0^*$ and $k = 1,2$, we have*

$$\int_{c_0} f^2(\mathbf{x})\,d\gamma_k(\mathbf{x}) = \sum_{u=1}^\infty \sum_{v=1}^\infty \Gamma(k,u,v) f(\mathbf{e}_u) f(\mathbf{e}_v).$$

PROOF. First, assume the $V_{n,j}$ are random. Then since (7.19) is verified, we also see for all $d \geq 1$ and $f \in c_{0,d}^*$ that

$$\lim_{n\to\infty} \sum_{i=1}^n E(f^4(\tilde{\mathbf{X}}_{n,i})) = 0.$$

Hence, by (7.5) and Lyapunov's central limit theorem, we have that $f(\tilde{\mathbf{S}}_{n,n})$ converges in distribution to a mean zero Gaussian random variable with variance given by (7.18) for all $d \geq 1$ and $f \in c_{0,d}^*$. If $\mu$ is a probability measure on the Borel subsets of $c_0$, and for all $k \geq 1, d \geq 1, f_1,\ldots,f_k \in c_{0,d}^*$, and $A$ is an arbitrary Borel set of $R^k$, then the probability distributions

$$F^{f_1,\ldots,f_k}(A) = \mu(\{\mathbf{x} \in c_0 : (f_1(\mathbf{x}),\ldots,f_k(\mathbf{x})) \in A\})$$

are the finite-dimensional distributions of $\mu$ on $c_0$ induced by $\bigcup_{d \geq 1} c_{0,d}^*$, and they uniquely determine $\mu$ on the Borel subsets of $c_0$. In view of the tightness obtained in Lemma 7.1, we thus have that $\tilde{\mathbf{S}}_{n,n}$ converges weakly to a unique probability on the Borel subsets of $c_0$, which for the moment we call $\mu$. What remains is to show that for every $f \in c_0^*$ this limiting measure makes $f$ a centered Gaussian random variable with variance determined by $\Gamma(1,\cdot,\cdot)$. Recalling that pointwise limits of centered Gaussian random variables are again centered Gaussian variables with limiting variances the limits of the variances, and that $\bigcup_{d\geq 1} c_{0,d}^*$ is weak star dense in $c_0^*$, it follows that $\mu$ is a centered Gaussian measure on $c_0$. Furthermore, if $f \in c_0^*$ and $f_d(\mathbf{x}) = f(\Pi_d(\mathbf{x})), \mathbf{x} \in c_0$, then for random $V_{n,j}$ we have

$$\int_{c_0} f^2(\mathbf{x})\,d\mu(\mathbf{x}) = \lim_{d\to\infty} \int_{c_0} f_d^2(\mathbf{x})\,d\mu(\mathbf{x}) = \lim_{d\to\infty} \sum_{u=1}^{d}\sum_{v=1}^{d} \Gamma(1,u,v)f(\mathbf{e}_u)f(\mathbf{e}_v).$$

Since $\sup_{i\geq 1} E(\xi_{n,i,j}^2) \leq c_{n,j}/k_{n,j}$, we have from (2.29), (2.26)–(2.28), and Cauchy–Schwarz that

$$\sup_{j_1,j_2\geq 1} |\Gamma(1,j_1,j_2)| < \infty.$$

Now $c_0^* = \ell_1$, and hence the dominated convergence theorem easily implies $\mu$ is a centered Gaussian measure on $c_0$ with covariance given by $\Gamma(1,\cdot,\cdot)$. Moreover, for each $f^* \in c_0$ we have

$$\int_{c_0} f^2(\mathbf{x})\,d\mu(\mathbf{x}) = \sum_{u=1}^{\infty}\sum_{v=1}^{\infty} \Gamma(1,u,v)f(\mathbf{e}_u)f(\mathbf{e}_v).$$

Hence, when $V_{n,j}$ is random, the centered Gaussian measure $\gamma_1$ exists as indicated, that is, its covariance is $\Gamma(1,\cdot,\cdot)$, and $\mu = \gamma_1$. Similarly, when the $V_{n,j}$ are replaced by $n$, then $\gamma_2$ exists as indicated, and $\mu = \gamma_2$. This last fact is easy to check by immediate simplifications of what we have done when $V_{n,j}$ is random, and the details are left to the reader. Hence for each choice of normalizers, there is a unique limiting Gaussian measure, and its finite-dimensional distributions are centered Gaussian measures determined by the appropriate covariance function. Therefore, the lemma is proved, and Theorem 2.4 is established. □

**Acknowledgments.** We would like to thank Mr. Bret Hanlon for help with simulations. Mr. Hanlon is a student of the second author. We also thank the referees and the associate editor for useful suggestions.

DEPARTMENT OF MATHEMATICS
UNIVERSITY OF WISCONSIN
MADISON, WISCONSIN 53706-1388
USA
E-MAIL: kuelbs@math.wisc.edu

DEPARTMENT OF STATISTICAL SCIENCE
CORNELL UNIVERSITY
ITHACA, NEW YORK 14853-4201
USA
E-MAIL: anv4@cornell.edu